\documentclass[10pt]{article}
\usepackage[english]{babel}

\usepackage{a4}
\usepackage{graphicx}
\usepackage{amsmath}
\usepackage{amsthm}
\usepackage{amssymb}

    \newtheorem{proposition}{Proposition}
    \newtheorem{theorem}{Theorem}
    \newtheorem{corollary}{Corollary}

   \theoremstyle{definition}
    \newtheorem{definition}{Definition}
 
    \newtheorem{example}{Example}
\theoremstyle{remark}
    \newtheorem{remark}{Remark}

\newcommand{\bb}{\mathbf}

\newcommand{\Pp}[1]{\bb{P}^{#1}}

\newcommand{\GL}{\mathrm{GL}}

\newcommand{\PGL}{\mathrm{PGL}}

\newcommand{\Z}{\bb{Z}}
\newcommand{\Q}{\bb{Q}}

\newcommand{\C}{\bb{C}}
\newcommand{\Fi}{\bb{F}}
\newcommand{\Ci}{\mathsf C}

\newcommand{\pu}{\bullet}

\DeclareMathOperator{\im}{Im}
\DeclareMathOperator{\Tr}{Tr}

\newcommand{\Pb}{\mathbf{P}}
\newcommand{\Ab}{\mathbf{A}}
\newcommand{\abs}[1]{\lvert#1\rvert}

\newcommand{\M}[1]{\mathcal M_{{#1}}}
\newcommand{\Mm}[2]{\mathcal M_{{#1},{#2}}}

\newcommand{\Hh}[1]{\mathcal H_{{#1}}}
\newcommand{\Hhm}[2]{\mathcal H_{{#1},{#2}}}

\newcommand{\Mb}[1]{\overline{\mathcal M}_{{#1}}}
\newcommand{\Mmb}[2]{\overline{\mathcal M}_{{#1},{#2}}}

\newcommand{\s}{\mathbb S}
\DeclareMathOperator{\Aut}{Aut}

\newcommand{\Eul}{\mathbf{e}}
\newcommand{\Eull}{\Eul_{\lambda}}
\newcommand{\etEul}{\Eul_{\acute{e}t}}
\newcommand{\etEull}{\Eul_{\acute{e}t, \lambda}}
\newcommand{\etEulu}{\Eul_{\acute{e}t, \mu}}
\newcommand{\Euls}[1]{\Eul^{\s_{#1}}}
\newcommand{\etEuls}[1]{\etEul^{\s_{#1}}}
\newcommand{\Ll}{\mathbf L}

\newcommand{\hh}[2][n]{{{}_{(#1)}}h_{#2}}
\newcommand{\V}{\mathcal V}

\newcommand{\mo}[3][\V]{{#1}(\!({#2},{#3})\!)}
\newcommand{\Mmo}[3][\MV]{{#1}(\!({#2},{#3})\!)}

\DeclareMathOperator{\Char}{\mathbb{C}\mathrm{h}}
\DeclareMathOperator{\ch}{ch}
\DeclareMathOperator{\Exp}{Exp}
\DeclareMathOperator{\Log}{Log}

\newcommand{\coh}[3][\Q]{H^{#2}({#3};{#1})}
\newcommand{\ut}[2]{u^{#1}t^{#2}}

\newcommand{\etcoh}[3][\Q_{l}]{H^{#2}_{\acute{e}t}({#3},{#1})}
\newcommand{\etcohl}[3][\Q_{l}]{H^{#2}_{\acute{e}t,\lambda}({#3},{#1})}

\DeclareMathOperator{\Fil}{Fil}

\newcommand{\ba}{\big|}
\newcommand{\X}{\mathcal{X}}
\DeclareMathOperator{\Conf}{Conf}
\newcommand{\B}[2]{B({#2},{#1})}
\newcommand{\tB}[2]{\tilde{B}({#2},{#1})}
\newcommand{\F}[2]{F({#2},{#1})}
\newcommand{\tF}[2]{\tilde{F}({#2},{#1})}

\newcommand{\op}{\mathring}
\newcommand{\duale}{\check{\ }}

\begin{document}

\title{The rational cohomology of $\Mb4$}

\author{Jonas Bergstr\"om\thanks{Institutionen f\"or matematik, Kungliga Tekniska h\"ogskolan, S--100 44 Stockholm, Sweden, E-mail address \texttt{jonasb@math.kth.se}} \and Orsola Tommasi\thanks{IMAPP, Radboud University Nij\-me\-gen, Toernooiveld 1, NL-6525 ED Nij\-me\-gen, the Netherlands, E-mail address \texttt{tommasi\@math.ru.nl}}  }

\date{November 9, 2006}

\maketitle
\begin{abstract}
We present two approaches to the study of the cohomology of moduli spaces of curves. Together, they allow us to compute the rational cohomology of the moduli space $\Mb4$ of stable complex curves of genus 4, with its Hodge structure.

\tableofcontents
\end{abstract}

\section{Introduction and results}

In this paper we compute the rational cohomology of the moduli space $\Mb4$ of stable complex curves of genus 4, together with its Hodge structure.

The moduli space $\Mb4$ has a natural stratification based on the topological type of the curves, or equivalently, on the dual graph of a curve. The locally closed substrata obtained in this way are finite quotients of products of moduli spaces of pointed non-singular curves. More precisely, if we denote by $\Mm gn$ the moduli space of non-singular curves of genus $g$ with $n$ marked points, we find that only the spaces $\Mm gn$ with $g\leq4$, $n\leq 8-2g$ are involved in the construction of the strata of $\Mb 4$. 

The moduli space $\Mb4$ is a proper smooth stack and satisfies therefore Poincar\'e duality. Hence, there are isomorphisms in all degrees between cohomology and cohomology with compact support, and all cohomology groups are pure with Hodge weight equal to the degree. For this reason, the rational cohomology of $\Mb4$ can be read off from the Hodge Euler characteristic of $\Mb4$, i.e., the Euler characteristic
$$\Eul(\Mb4):=\sum_{j \in \Z}(-1)^j[H_c^j(\Mb4;\Q)]$$ 
of $H_c^\pu(\Mb4;\Q)$ in the Grothendieck group of rational mixed Hodge structures. 

Since the Hodge Euler characteristic is additive, it is possible a priori to determine the Hodge Euler characteristic of $\Mb4$ from that of its strata. As was said above, each stratum is the quotient of a product of spaces of the form $\Mm gn$ by the action of a finite group. The Hodge Euler characteristic of all $\Mm gn$'s occurring in the strata does not give enough information to compute the Hodge Euler characteristic of the strata, since it does not keep track of the action of the group. What we need instead is the $\s_n$-equivariant Hodge characteristic of the $\Mm gn$'s, which also encodes the structure of the cohomology groups as representations of the symmetric group $\s_n$, which acts naturally on $\Mm gn$ by permuting the marked points. 

Once all $\s_n$-equivariant Hodge Euler characteristics of the spaces $\Mm gn$ are known for $g\leq4$, $n\leq 8-2g$, determining the Hodge Euler characteristic of $\Mb 4$ becomes a problem in the combinatorics of graphs. Indeed, this is the approach of Getzler and Kapranov in \cite{GK}, where they develop the theory of \emph{modular operads} and apply it to the homology of the moduli spaces of curves. In particular, this enables them to give a formula (\cite[8.13]{GK}), that expresses the relationship of Euler characteristics between $\Mmb gn$ and the $\Mm {\tilde g}{\tilde n}$ appearing in its stratification. As remarked in \cite{G-semi}, the formula applies to Hodge Euler characteristics, and in general to all Euler characteristics taking values in the Grothendieck group of a symmetric monoidal category that is additive over a field of characteristic 0 and has finite colimits. 

The computation of the cohomology of $\Mb 4$ is thus essentially reduced to the computation of the equivariant Hodge Euler characteristic of all spaces $\Mm gn$ with $g\leq4$, $n\leq 8-2g$. The equivariant Hodge characteristics of all moduli spaces $\Mm 0n$ and $\Mm 1n$, and of $\Mm 2n$ ($n\leq3$), $\Mm 3n$ ($n\leq1$) are already known (see \cite{G-0}, \cite{G-res}, \cite{G-2}, \cite{Looij}, \cite{GL}). In this paper, we present two methods that allow us to calculate the Hodge Euler characteristic in the remaining cases. As a result, we establish the following.

\begin{theorem}\label{mainthm}
The Hodge Euler characteristic of $\Mb 4$ is equal to
$$\Ll^9+4\Ll^8+13\Ll^7+32\Ll^6+50\Ll^5+50\Ll^4+32\Ll^3+13\Ll^2+4\Ll+\mathbf{1}$$ 
where $\Ll$ denotes the class of the Tate Hodge structure of weight $2$ in the Grothen\-dieck group of rational Hodge structures. 
\end{theorem}

We notice that $\Mb 4$ has no odd cohomology. This is in agreement with the result in the article \cite{Arb-Cor} by Arbarello-Cornalba, that the first, third and fifth cohomology group of $\Mmb gn$ is zero for all $g$ and $n$. Furthermore, the theorem suggests the possibility that all cohomology classes of $\Mb4$ are tautological. This has indeed been proved recently by Carel Faber and Rahul Pandharipande (\cite{tautological}).

The first method we present is based on the equivariant count of the number of points of $\Mm gn$ defined over finite fields. It has previously been used to get cohomological information on the moduli space of curves of genus two by Faber-Van der Geer in \cite{Faber-Geer1}, \cite{Faber-Geer2} and on the moduli spaces of pointed curves of genus 0 by Kisin-Lehrer in \cite{Lehrer}.

Here we use it to get the wanted information for curves of genus two and three. The key remark is that if we consider the number of points of $\Mm 2n$ or $\Mm 3n$ defined over a finite field $k$ as a function of the number $q$ of elements of $k$, the function we get is a polynomial in $q$ if $n$ is sufficiently small. 
This is not merely true for the number of points of $\Mm gn$, i.e., the trace of Frobenius on the \'etale cohomology with compact support of $(\Mm gn)_{\overline{k}}$. Also the trace of the composition of Frobenius with any automorphism of \'etale cohomology induced by the action of an element of $\s_n$ on $\Mm gn$ is a polynomial in $q$. We say that in these cases, the \emph{equivariant count} of the number of points of $\Mm gn$ gives a polynomial. This is true for all moduli spaces occurring in the boundary of $\Mb4$. The equivariant count of points of $\Mm gn$ for $g\in\{2,3\}$ is performed for a wider range of indexes $n$ and explained in greater detail by the first author in \cite{Jonas1} and \cite{Jonas2}.

The equivariant count of the number of points of $\Mm 2n$ is based on the fact that all curves of genus two are hyperelliptic. In general, if the characteristic of $k$ is odd, the moduli space $\Hh g$ of hyperelliptic curves of a fixed genus $g\geq2$ is isomorphic to the moduli space of binary forms of degree $2g+2$ without multiple roots. This allows us to reduce the equivariant count of the number of points defined over a finite field $k$ of the moduli space of genus $g$ hyperelliptic curves with $n$ marked points, to a computation involving only the quadratic characters of finite extensions of $k$ of the values of square-free monic polynomials of the appropriate degrees, at points on the projective line.

The equivariant count of the number of points of $\mathcal Q_{n}:=\Mm 3n\setminus\Hhm 3n$ is based on the fact that this space coincides with the moduli space of non-singular quartic curves in $\Pp2$. Hence, in order to know the number of points defined over a finite field $k$ of $\mathcal Q_{n}$ it is enough to count the number of $(n+1)$-tuples $(C,p_1,\dots,p_n)\in\Pp{14}(k)\times({\Pp2(k)})^n$, where $\Pp{14}(k)$ is identified with the space of quartic curves in $\Pp2$ with coefficients in $k$, $C$ is a non-singular curve and $p_1,\dots,p_n$ are distinct points on $C$. It turns out that it is easier to count the number of singular curves than the number of non-singular ones.

Hence, it is possible to obtain the equivariant count of points in $\mathcal Q_{n}$ by counting, for each choice of an $n$-tuple defined over $k$ of distinct points in $\Pp2$, the number of singular curves passing through these points. The latter count is performed by a modified version of the sieve principle.

From the fact that the equivariant count of points of $\Mm gn$ gives polynomials in the cases we investigate, it follows that we can use these polynomials to obtain the Hodge Euler characteristic of the moduli spaces. A theorem by Van den Bogaart and Edixhoven (\cite{EB}) ensures that the cohomology of $\Mmb gn$ is pure and of Tate Hodge type whenever the number of points of $\Mmb gn$ is a polynomial in $q$. Furthermore, the Betti numbers and the Hodge weight of the cohomology classes are given by the coefficients of the polynomial. It is not difficult to extend this result so to keep track of the action of $\s_n$ on the cohomology of $\Mmb gn$. Note that if the equivariant count of points gives polynomials for all $\Mm{\tilde g}{\tilde n}$ occurring in the stratification of $\Mmb gn$, also the count of the number of points of $\Mmb gn$ will give a polynomial. This implies that in this situation we can determine the $\s_n$-equivariant Hodge characteristic of $\Mm gn$ from the $\s_n$-equivariant count of the number of its points. Specifically, this allows us to establish the following.

\begin{theorem}\label{char24}
The equivariant Hodge Euler characteristic (written in terms of the Schur polynomials, see Section \ref{PartI}) of $\Mm 24$ is equal to $$(\Ll^7+\Ll^6-\Ll^5-\Ll^4-\Ll^2+\Ll+\mathbf{1})  s_4 + (\Ll^6-\Ll^5-\Ll^2)  s_{31}- \Ll^2  s_{2^2} + (\Ll^2-\mathbf{1})  s_{21^2}.$$
\end{theorem}

\begin{theorem}\label{char32}
The equivariant Hodge Euler characteristic of $\Mm 32$ is equal to $$(\Ll^8+2\Ll^7+\Ll^6-\Ll^5+\Ll^2+\Ll) s_2 + (\Ll^7+\Ll^6+\Ll) s_{1^2}.$$
\end{theorem}

The second method we present is based on the fact that $\M 4$ and $\M 3$ have a stratification such that all strata can be interpreted as quotients of complements of discriminants in a complex vector space by the action of a reductive group. On the one hand, the rational cohomology of the complement of such a discriminant can be computed by Vassiliev-Gorinov's method (see \cite{Vart}, \cite{Gorinov}). On the other hand, by a generalization by Peters and Steenbrink of the Leray-Hirsch theorem (\cite{PS}), the rational cohomology of these complements of discriminants is the tensor product of the rational cohomology of the quotient and that of the group acting. As a consequence, if the cohomology of the complement of the discriminant is known, it is straightforward to obtain the cohomology of the quotient. We recall here the version of Vassiliev-Gorinov's method used in \cite{OT} to compute the rational cohomology of $\M 4$. Moreover, we show that it is easy to apply Peters-Steenbrink's generalized Leray-Hirsch theorem to incidence correspondences. This implies that we can use the same methods to compute the rational cohomology of moduli spaces of non-singular pointed curves. As a first easy example, we present here the calculation of the rational cohomology of $\Mm 31$. The rational cohomology of $\Mm 31$ is given in Corollary~\ref{M31}, which refines the result given in \cite{GL}, which corrects an error in \cite[4.10]{Looij}.

The outline of the paper is as follows. In Section~\ref{PartI}, we present the formula of Getzler and Kapranov on the generating function of the Euler characteristic of modular operads, and use it to establish Theorem~\ref{mainthm} from the Hodge Euler characteristics of the spaces $\Mm gn$ with $g\leq 4$, $n\leq 8-2g$. In Section~\ref{PartIV}, following~\cite{EB}, we explain why the equivariant count of points of $\Mm gn$ for $g\leq 3$, $n\leq 8-2g$ gives the equivariant Hodge Euler characteristic of these spaces. This requires us to use an equivariant version of \cite[Theorem 2.1]{EB}. In Section~\ref{PartIII}, we introduce the methods used for the equivariant counts of the number of points of $\Mm gn$ in the cases we are interested in. In Section~\ref{PartII} we present the relation between the study of the cohomology of moduli spaces of curves of genus 3 and of genus 4 and the cohomology of the complement of certain discriminants. Moreover we introduce the generalized Leray-Hirsch theorem and explain why it is applicable also in the case of moduli spaces of curves with marked points. We recall the result in \cite{OT} on the rational cohomology of $\M4$, and compute the rational cohomology of $\Mm 31$. These computations are achieved using Vassiliev-Gorinov's method, which is presented in Section~\ref{VG}.

\subsubsection*{Acknowledgments}
The cooperation between the two authors started with a visit of the second author at KTH in Fall 2004. The second author would like to thank Carel Faber for the invitation and KTH for the hospitality. We would like to thank Joseph Steenbrink and Carel Faber for useful comments. 
Moreover, we would like to thank Theo van den Bogaart and Bas Edixhoven for allowing us to read early versions of their article \cite{EB}.

\subsection*{Notation}
Let $\Mm{g}{n}$ be the moduli space of irreducible, non-singular, projective curves of genus $g$ with $n$ distinct marked points. Furthermore, let $\Mmb{g}{n}$ be the moduli space of stable curves of genus $g$ with $n$ distinct marked points. Both these moduli spaces are smooth Deligne-Mumford stacks defined over $\Z$. The symmetric group $\s_n$ of permutations of $n$ elements acts on both $\Mm gn$ and $\Mmb gn$ by permuting the $n$ marked points on the curves. 

If $\X$ is a Deligne-Mumford stack and $S$ is a scheme, we denote by $[\mathcal X(S)]$ the set of isomorphism classes of the category $\mathcal X(S)$.

For every pair of non-negative integers $g,n$ such that $2g+n-2>0$, we let $$\mathcal D_{g,n}:=\{(\hat g,\hat n): 0\leq \hat g\leq g, \max\{0,3-2\hat g\}\leq \hat n\leq 2(g-\hat g) + n\}.$$

To state the results in this paper, we need to consider the Grothendieck groups of several categories. We will denote by $\mathsf{MHS_\Q}$ the category of rational mixed Hodge structures, and by $\mathsf{Gal}$ the category of $\mathrm{Gal}(\bar\Q/\Q)$-representations. For any abelian category $\mathsf C$, we will denote by $\mathrm K_0(\mathsf C)$ the Grothendieck group of $\mathsf C$. Note that $\mathrm K_0(\mathsf{MHS}_\Q)$ coincides with $\mathrm K_0(\mathsf{HS}_\Q)$, the Grothendieck group of (pure) Hodge structures. 

We denote by $\Q(-k)$ the one-dimensional rational Tate Hodge structure of weight $2k$. The class of $\Q(-1)$ in $\mathrm K_0(\mathsf{HS}_\Q)$ is denoted by $\Ll$. For any rational mixed Hodge structure $H$, we denote by $H(-k)$ the Tate twist $H\otimes\Q(-k)$. Similarly, we denote by $\Q_l(k)$ the $k$-th Tate twist of the trivial $\mathrm{Gal}(\bar\Q/\Q)$-representation.

\section{The stratification of $\Mmb GN$}\label{PartI}
The strategy of this paper is to compute the rational cohomology of the moduli space $\Mb 4$ of stable curves (over the complex numbers), from what we know on the cohomology of some moduli spaces of smooth $n$-pointed curves of genus $g\leq 4$.

The stratification based on the dual graph gives a way to divide every moduli space $\Mmb GN$ into strata which are explicitly related to the moduli spaces $\Mm gn$ with $(g,n)\in\mathcal D_{G,N}$. This gives the intuition that all information on $\Mmb GN$ can be deduced from the appropriate information on moduli spaces of smooth pointed curves.
 
The aim of this section is to present in Theorem~\ref{thmGK} a formula by Getzler and Kapranov (\cite{GK}) which expresses the relationship between the $\s_N$-equivariant Euler characteristic of $\Mmb GN$ and the $\s_n$-equivariant Euler characteristics of the $\Mm gn$'s with $(g,n)\in\mathcal D_{G,N}$. Then we will apply this theorem to $\Mb 4$, in the case of Hodge Euler characteristics, to give a proof of Theorem \ref{mainthm}.

First, we explain what the dual graph of a stable curve is. Let $(C,p_1,\dots,p_N)$ be a stable curve of genus $G$ with $N$ marked points. 
Its dual graph is the labelled graph $\Gamma$ such that:\begin{itemize}
\item Every vertex $v_i$ of $\Gamma$ corresponds to an irreducible component $C_i$ of $C$, and its label $g_i$ is the geometric genus of $C_i$.
\item The edges of $\Gamma$ correspond to the nodes of $C$ (e.g., two distinct vertices are joined by an edge if and only if the corresponding components intersect).
\item There are $N$ half-edges, labelled from $1$ to $N$, corresponding to the $N$ marked points.
\end{itemize}
For instance, a non-singular $N$-pointed curve of genus $G$ corresponds to a tree with one vertex of genus $G$ and $N$ half-edges. 

Finally, an automorphism of a labelled graph is defined to be an automorphism of the underlying non-labelled graph that preserves the labelling of the vertices and that fixes each half-edge.

For any connected labelled graph $\Gamma$ with $N$ half-edges and such that 
$$\sum_i(g_i-1)+\#\{Edges\}+1=G,$$
we can consider the moduli space $\M{}(\Gamma)$ of curves with dual graph $\Gamma$. Each $\M{}(\Gamma)$ is locally closed in $\Mmb GN$. Note that the spaces $\M{}(\Gamma)$ can be expressed in terms of moduli spaces of non-singular curves. Indeed, $\M{}(\Gamma)\cong\prod_i{\Mm{g_i}{n_i}}/\Aut(\Gamma)$, where $g_i$ denotes, as usual, the label of the $i$-th vertex $v_i$ of $\Gamma$, and $n_i$ is the sum of the number of edges and half-edges incident to the vertex $v_i$. 

Since $\Mmb GN$ is complete and satisfies Poincar\'e duality, the knowledge of the cohomology of $\Mmb GN$ as a graded vector space with mixed Hodge structure is equivalent to that of the \emph{Hodge Euler characteristic} of $\Mmb GN$. That is, the Euler characteristic of $H_c^\pu(\Mmb GN;\Q)$ in the Grothendieck group of rational mixed Hodge structures,
$$\Eul(\Mmb GN)
:=\sum_{j\in\Z}(-1)^j[H_c^j(\Mmb GN;\Q)]\in 
\mathrm K_0(\mathsf{MHS}_\Q).$$

This Euler characteristic does not contain all the information we need.
Namely, it does not give information on the action of the symmetric group $\s_n$, which acts on $\Mm gn$ (and $\Mmb gn$) by permuting the marked points. In order to express the induced action on cohomology, we will make use of symmetric functions in our notation, given the well-known correspondence between them and the characters of the symmetric group.

We denote by $\Lambda_n$ the ring of symmetric functions in $n$ variables. The ring $\Lambda_n$ is generated by the complete symmetric functions $\hh1,\dots,\hh n$, where $\hh j:=\sum_{1\leq i_1\leq\dots\leq i_j\leq n} x_{i_1}x_{i_2}\dots x_{i_j}$.
There is a natural projection $\Lambda_{n+1}\rightarrow \Lambda_n$ that maps $\hh[n+1]j$ to $\hh j$ if $j\leq n$ and $\hh[n+1]{n+1}$ to $0$. This allows us to define 
$$\Lambda:={\underleftarrow{\lim}}\ \Lambda_n.$$

Note that the elements of $\Lambda$ are infinite series in the variables $\{x_j\}_{j\geq1}$, and can be written as (possibly infinite) sums of \emph{complete symmetric functions}
$h_j:=\sum_{1\leq i_1\leq\dots\leq i_j}x_{i_1}\dots x_{i_j}$. For every $n$, the ring $\Lambda_n$ can be identified with the subring $\Z[h_1,\dots,h_n]$ of $\Lambda$.

The correspondence between elements of $\Lambda_n$ and representations of $\s_n$ (for some fixed $n$) is given by the Schur polynomials. The Schur polynomials of degree $n$ are in one to one correspondence with partitions of $n$ and we denote by $s_\lambda$ the Schur polynomial corresponding to the partition $\lambda=(\lambda_1,\dots,\lambda_s)$ (with $\lambda_1\geq\dots\geq\lambda_s\geq1$). As is well known, the set of Schur polynomials freely generates $\Lambda$ as an abelian group.

We will also consider another set of symmetric functions, namely the \emph{power sums} $p_n:=\sum_i x_i^{n}$, which constitute a $\Q$-basis of $\Lambda\otimes\Q$. 

In order to compute $\Eul(\Mb4)$, we use the formula by Getzler and Kapranov given in Theorem~\ref{thmGK}. Its formulation and proof are based on the theory of \emph{modular operads}. Since there is an excellent exposition of it in \cite{GK}, we will not explain this theory here, but refer to \cite{GK} for all definitions and properties we use. Following \cite{G-res}, we will formulate the theory for any symmetric monoidal category $\mathsf C$ that is additive over a field of characteristic 0 and has finite colimits. For the definitions of these concepts see for instance \cite[Chapter VII, Chapter XI]{tensor} and \cite{G96}. Note that $\mathsf{MHS}_\Q$ is an example of such a category.

To any such category $\mathsf C$ we associate the category $\mathsf C^{\s_n}$ whose objects are objects of $\mathsf C$ equipped with an action of the symmetric group $\s_n$ by morphisms in $\mathsf C$. 

\begin{theorem}[{\cite[Theorem 4.8]{G96}}]
There is a canonical isomorphism $$\mathrm K_0(\Ci^{\s_n})\cong\mathrm K_0(\Ci)\otimes\Lambda_n.$$
\end{theorem}

\begin{definition}[\cite{GK}]
A stable $\s$-module $\V$ in the category $\mathsf C$ is a collection, for all $g,n\geq 0$, of chain complexes $\{ \mo gn_i \}$ of objects of $\mathsf C^{\s_n}$ such that  $\mo gn=0$ if $2g+n-2\leq0$. 
\end{definition}

\begin{definition}\label{d:chn}\label{eqHC}
Let $R:=\{R_i\}$ be a finite chain complex of objects of $\mathsf C^{\s_n}$ for some $n\geq0$. The characteristic of $R$ is defined as 
$$\ch_n(R):= \sum_i (-1)^i [R_i] \in \mathrm K_0(\mathsf C^{\s_n}) \cong \mathrm K_0(\mathsf C) \otimes \Lambda_n.$$ 

When $\mathcal X$ is an algebraic stack over $\C$, with an action of the symmetric group $\s_n$, the rational cohomology with compact support of $\mathcal X$ has a natural structure of a chain complex $\{C_i\}$ of objects of $\mathsf{MHS}_\Q^{\s_n}$ by setting $C_i:=H_c^i(\mathcal X;\Q)$ and defining all differentials to be zero. 
We define the \emph{$\s_n$-equivariant Hodge Euler characteristic} $\mathcal X$ as
$$\Euls n(\mathcal X):=\ch_n(H_c^\pu(\X;\Q))\in\mathrm{K}_0(\mathsf{MHS}_\Q)\otimes\Lambda_n.$$ 
\end{definition}

\begin{definition}\label{d:Char}
Let $\V$ be a stable $\s$-module. Then the \emph{characteristic of $\V$} is defined as 
$$\Char(\V):=\sum_{2g+n-2>0}\hbar^{g-1}\ch_n(\mo gn)\in \mathrm K_0(\mathsf C) \otimes \Lambda(\!(\hbar)\!),$$
where $\Lambda(\!(\hbar)\!)$ is the ring of Laurent series with coefficients in $\Lambda$.
\end{definition}

\begin{remark} On $\Lambda$ there is a unique associative operation, called \emph{plethysm}, satisfying the conditions\begin{enumerate}
\item $(f_1+f_2)\circ g = f_1\circ g+f_2\circ g$;
\item $(f_1f_2)\circ g=(f_1\circ g)(f_2\circ g)$;
\item If $f=f(p_1,p_2,\dots)$ then $p_n\circ f=f(p_n,p_{2n},\dots)$.
\end{enumerate}

Following \cite{G96}, the plethysm $\circ$ can be extended to a map $\Lambda \times (\mathrm{K}_0(\Ci)\otimes\Lambda)\rightarrow \mathrm{K}_0(\Ci)\otimes\Lambda$ by setting $h_n\circ [R]$ for every $n\geq 0 $, $R\in\Ci$,  to be the class in $\mathrm K_0(\Ci)$ of the $n$-th symmetric product of $R$. In particular, when $\Ci$ is the category of rational mixed Hodge structures, we have $p_n\circ \Ll = \Ll^n$, where $\Ll$ is the class of the Tate Hodge structure of weight $2$. We also extend the plethysm to $\Lambda(\!(\hbar)\!)$ by letting $p_n\circ\hbar=\hbar^n$.
\end{remark}

\begin{definition}\label{d:expetc}\begin{enumerate}
\item
Suppose $f=\sum_{\alpha\in\Z}f_\alpha[R_{\alpha}]\hbar^\alpha\in\mathrm{K}_0(\Ci)\otimes\Lambda(\!(\hbar)\!)$, where $f_\alpha \in \Lambda$ and $[R_{\alpha}] \in\mathrm{K}_0(\Ci) $, is such that for all $\alpha\in\Z$ all monomials occurring in $f_\alpha$ have degree at least $1-2\alpha$.

Then the plethystic exponential of $f$ is defined by the expression:
$$\Exp(f)=\sum_{n=0}^\infty h_n\circ f.$$
\item We denote the inverse of the plethystic exponential by:
$$\Log(f)=\sum_{n=1}^\infty\frac{\mu(n)\log(p_n)\circ f}n.$$
\item On $\Lambda(\!(\hbar)\!)$, an analogue of the Laplacian is given by: 
$$\Delta=\sum_{n=1}^\infty \hbar^n \left(\frac n2 \frac{\partial^2}{\partial 
p_n^2}+\frac\partial{\partial p_{2n}}\right).$$
\end{enumerate}
\end{definition}

Recall from \cite{GK} that it is possible to associate with every stable $\s$-module the free modular operad $\mathbb M\V$ generated by $\V$. The following theorem gives the relationship between the characteristics of the stable $\s$-modules $\V$ and $\mathbb M\V$.

\begin{theorem}[{\cite[Thm 8.13]{GK}}]\label{thmGK}
Let $\V$ be an $\s$-module in the category $\Ci$. Then
$$
\Char (\mathbb M\V)=\Log(\exp(\Delta)\Exp(\Char(\V)).
$$
\end{theorem}

\begin{proof}[of Theorem~\ref{mainthm}]
Let us consider $\mo gn := H_c^\pu (\Mm gn;\Q)$. 
From the results in \cite{GK} it follows that $\V$ is a modular operad in the category of stable $\s$-modules with mixed Hodge structures. The free modular operad associated with $\V$ is $H_c^\pu(\Mmb gn;\Q)$. In this situation, Theorem~\ref{thmGK} can be used to compute the generating function 
$$\Char(\mathbb M\V)=\sum_{2g+n-2>0}\hbar^{g-1}\Euls n(H_c^\pu(\Mmb gn;\Q)),$$
from the generating function $\Char(\V)$.

The Hodge Euler characteristic of $\Mb 4$ is exactly the part of the coefficient of $\hbar^3$ in $\Char(\mathbb M\V)$ with degree 0 in $\Lambda$. It is not necessary to compute the whole $\Char (\mathbb M\V)$ in order to obtain this. Since the Euler characteristic $2-2g-n$ of the dual graph behaves additively throughout the whole computation, only the $\mo gn$'s with Euler number not less than $-6=2-2\cdot 4$ give a contribution. Therefore, all we need to know is the $\s_n$-equivariant Hodge Euler characteristic of $\Mm gn$ for all $(g,n)$ in $\mathcal D_{4,0}$. 

For $g=0$, $g=1$, $g=2$ and $n\leq3$, and $g=3$, $n\leq1$, these results are already known (see \cite{G-0}, \cite{G-res}, \cite{G-2}, \cite{Looij}, \cite{GL}). 
The equivariant Hodge Euler characteristic of $\Mm 24$ and $\Mm 32$ are given in Theorems~\ref{char24} and~\ref{char32} respectively. These theorems follow from the results in Sections~\ref{PartIV} and~\ref{PartIII} of this paper and rely on the equivariant count of points made by the first author in \cite{Jonas1} and \cite{Jonas2}.

The rational cohomology of $\M 4$ has been computed by the second author in \cite{OT} with the methods explained in Sections~\ref{PartII} and \ref{VG}. In particular, Theorem~\ref{cohM4} implies that $\Eul(\M4)=\Ll^9+\Ll^8+\Ll^7-\Ll^6$.

Note that for all $(g,n)\in\mathcal D_{4,0}$, the $\s_n$-equivariant Hodge Euler characteristic of $\Mm gn$ lies in $\Q[\Ll]\otimes\Lambda$. Hence, our computation only involves polynomials in $\Ll$ and certain symmetric functions. 

To obtain the result, we used a computer program (by Del Ba\~{n}o), based on Theorem~\ref{thmGK}, which computes the part of $\Char(\mathbb M\V)$ of degrees $G-1$ in $\hbar$ and $N$ in $\Lambda$, whenever $\ch_n(\mo {g}{n})$ is given and is a polynomial in $\Ll$ for every $(g,n)$ in $\mathcal D_{G,N}$. This computer program makes extensive use of Stembridge's package \textsf{SF} (\cite{SF}) for computations with symmetric functions.
\end{proof}

\section{From counting to cohomology}\label{PartIV}
For every integer $q$ which is the power of a prime number, denote by $\Fi_q$ the finite field with $q$ elements and by $F_q$ its Frobenius map. Denote by $(M_{g,n})_{\bar\Fi_q}$ the coarse moduli space of the stack $(\Mm gn)_{\bar\Fi_q}$. Fix two non-negative integers $G$ and $N$ such that $2G+N-2>0$, and assume that the functions
$$f_{g,n,\sigma}(q)= \abs{(M_{g,n})_{\bar\Fi_q}^{F_q \cdot \sigma}}$$
are polynomials in $q$ for every $(g,n)$ in $\mathcal D_{G,N}$ and for every $\sigma\in\s_n$.

Then we will see that we can determine the {\'e}tale cohomology of $(\Mmb GN)_{\overline{\Q}}$ with its structure as $\mathrm{Gal}(\bar\Q/\Q)$- and $\s_N$-representation. From this we can determine the action of the symmetric group and the Hodge structure on the Betti cohomology of the complex stack $(\Mmb GN)_{\C}$. 

\subsection{{\'E}tale cohomology} \label{etGK}
In this section we will apply the theory of Getzler and Kapranov, presented in Section \ref{PartI}, to the stacks $\Mm gn$ in the setting of $l$-adic {\'e}tale cohomology.

Let $k=\Fi_q$ denote the finite field with $q$ elements. All varieties and stacks mentioned in the following will be assumed to be defined over $\overline{k}$ if nothing else is specified. 

We denote by $\etcoh{\bullet}{-}$ compactly supported $l$-adic {\'e}tale cohomology. Observe that the $c$, standing for compactly supported, is omitted in this notation. Unless otherwise stated we assume that $l$ does not divide $\abs{k} = q$. Let $\mathsf{Gal}_{\Q_l}$ be the category of $\Q_l$-vector spaces equipped with the $l$-adic topology that have a continuous action of the absolute Galois group $\mathrm{Gal}(\bar \Q / \Q)$. This is a symmetric monoidal additive category over $\Q_l$, with finite colimits. Hence, Theorem \ref{thmGK} holds in this category. 

\begin{definition}
Suppose we have a space $\X$ with an action of $\s_n$ commuting with the action of the absolute Galois group. Denote by $\etcohl{j}{\X}$ the Galois subrepresentation of $\etcoh{j}{\X}$ which is the sum of all copies of the irreducible representation of $\s_n$ indexed by the partition $\lambda$ that appear in $\etcoh j\X$. Then for every partition $\lambda$ of $n$ we define $\etEull(\X)$ to be the Euler characteristic of $\etcohl{j}{\X}$ in the Grothendieck group of $\mathsf{Gal}_{\Q_l}$.
\end{definition}

Define a stable $\s$-module $\V$ in $\mathsf{Gal}_{\Q_l}$ by  
$$\mo gn := \etcoh{\bullet}{\Mm gn}.$$ 

Taking the characteristic of $\mo gn$ in the Grothendieck group of equivariant Galois representations gives
$$\ch_n \bigl( \mo gn \bigr) = \sum_{\lambda \,
  \vdash n} (\chi_\lambda(id))^{-1}\etEull (\Mm gn) s_{\lambda} \in
  \mathrm{K}_0(\mathsf{Gal}_{\Q_l}) \otimes \Lambda_n .$$ 
The properties of Euler characteristics now ensure that 
$$\ch_N \bigl( \Mmo GN \bigr) = \sum_{\lambda \,
  \vdash N} (\chi_\lambda(id))^{-1}\etEull (\Mmb GN) s_{\lambda} 
\in \mathrm{K}_0(\mathsf{Gal}_{\Q_l}) \otimes \Lambda_N.$$ 

Applying Theorem \ref{thmGK} to $\V$, we can express $\etEull(\Mmb GN)$ in terms of $\etEulu(\Mm gn)$ for partitions $\mu$ of $n$ and for $(g,n) \in \mathcal{D}_{G,N}$. Note that the plethysm of $p_j\in\Lambda$ with a Galois representation of the form $\Q_l(i)$ for some $i \in \Z$ equals $[\Q_l(i\cdot j)]\in\mathrm K_0(\mathsf{Gal}_{\Q_l})$. Hence, if $\etEulu(\Mm gn)$ is a sum of Tate twists of the trivial Galois representation for all partitions $\mu$ of $n$ and for all $(g,n) \in \mathcal{D}_{G,N}$ then this also holds for $\etEull(\Mmb GN)$. 

\subsection{The Lefschetz trace formula} \label{count-to-trace}
Let $F_q$ be the geometric Frobenius map belonging to the finite field $k$. Let $\X_k$ be a smooth Deligne-Mumford stack of constant dimension and of finite type over $k$ that has an action of $\s_n$. Denote by $X_{\bar k}$ the coarse moduli space of $\X_{\bar k}:=\X_k \otimes_k \bar{k}$. 

By an $\s_n$-\emph{equivariant count} of the number of points defined over $k$ of $\X_{\bar k}$ we mean a count, for each $\sigma \in \s_n$, of the number of fixed points of $F_q \cdot \sigma$ acting on $X_{\bar k}$. Note in particular that these numbers are the same for all $\sigma$ with the same cycle type. 

The Lefschetz trace formula, generalized in \cite[3.1.2]{Behrend} to Deligne-Mumford stacks (see also \cite[Section~2]{Behrthesis}), gives the equality
$$\mathrm{Tr}\bigl( F_q \cdot \sigma,
\etEul (\X_{\bar k}) \bigr) = \abs {X_{\bar k}^{F_q\cdot \sigma}}.$$ 

Let $\chi_{\lambda}$ be the character of the irreducible representation of $\s_n$ indexed by $\lambda$. Then the endomorphism
$$\pi_\lambda:=\frac{1}{n!}  \chi_{\lambda}(id)  \sum_{\sigma \in \s_n}
\chi_{\lambda}(\sigma)  \sigma$$ 
is the projection of $\etcoh i {\X_{\bar k}}$ onto $\etcohl i{\X_{\bar k}}$ (see for instance \cite[2.31]{Fulton-Harris}). Therefore,
\begin{equation} \label{eq-trace}
\frac{1}{n!} \chi_{\lambda}(id) \sum_{\sigma \in \s_n} \chi_{\lambda}(\sigma) \abs{X_{\bar k}^{\sigma \cdot F_q }} =  
\mathrm{Tr}\bigl( F_q \cdot \pi_\lambda,\etEul(\X_{\bar k}) \bigr) =\mathrm{Tr} \bigl( F_q, \etEull (\X_{\bar k}) \bigr) \end{equation} 
which gives a formula expressing $\Tr(F_q,\etEull(\X_{\bar k}))$ as a function of the $\abs{X_{\bar k}^{F_q \cdot \sigma}}$ with $\sigma\in\s_n$.

\subsection{The Galois action} \label{count-to-adic}
In this section we will present an equivariant version of Theorem 2.1 in \cite{EB}. 

\begin{theorem} \label{thm-count-adic}
Let $\mathcal{\X}$ be a Deligne-Mumford stack defined over $\Z$ which is proper, smooth, of pure relative dimension $d$ and that has an action of $\s_n$. Let $X_{\bar{\Fi}_p}$ be the coarse moduli space of $\mathcal{\X}_{\bar{\Fi}_p}$.

For every partition $\lambda$ of $n$, denote by $\chi_{\lambda}$ the character of the irreducible representation of $\s_n$ indexed by $\lambda$. Furthermore let $S$ be a set of primes of Dirichlet density 1. 

Assume that for a partition $\lambda$ of $n$ there exists a polynomial 
$P_\lambda(t) \in \Q[t]$ such that
\begin{equation} \label{cond-adic}  
\frac{1}{n!}  \chi_{\lambda}(id)   \sum_{\sigma \in \s_n} 
\chi_{\lambda}(\sigma) \abs{X_{\bar{\Fi}_p}^{\sigma \cdot F_{p^r}}} = P_\lambda(p^r)
\end{equation}
for all $r \in \Z_{\geq 1}$ and $p \in S$. 

Then $P_\lambda(t)$ has degree $d$ and non-negative integer coefficients. Furthermore, if we let $b_j$ be the coefficient of $q^j$ in $P_\lambda$, then for all primes $l$ and all $i \geq 0$ there is an isomorphism of $\mathrm{Gal}(\bar{\Q}/\Q)$-representations 
$$\etcohl{i}{\X_{\bar{\Q}}} \simeq \begin{cases} 0 & \text{if $i$ is odd} 
\\ \Q_l (-i/2)^{b_{i/2}} & \text{if $i$ is even} \end{cases}$$ 
\end{theorem}

\proof
The proof follows very closely that of Theo van den Bogaart and Bas Edixhoven in \cite{EB}, thus we will only include the vital steps. 

Assume Condition~\eqref{cond-adic} holds for a partition $\lambda$. 
Using equation (\ref{eq-trace}) we get, for some $b_{i} \in \Q$,
\begin{equation} \label{eq-poly}
\mathrm{Tr} \bigl( F_{p^r}, \etEull (\X_{\bar{\Fi}_p}) \bigr) = \sum_i b_{i} p^{ri}
\end{equation} 
for all $r \in \Z_{\geq 1}$ and $p \in S$. Note here that $\mathrm{Tr}(F_{p^r},[\Q_l(-i)]) = p^{ri}$.

In the proof of Theorem~2.1 of~\cite{EB} it is shown that (if $p\neq l$) the Galois representation $\etcoh i {\X_{\bar{\Fi}_p}}$ is unramified and all eigenvalues of the Frobenius map $F_{p^r}$ have complex absolute value $p^{ri/2}$. Since $\etcohl{i}{\X_{\bar{\Fi}_p}}$ is a Galois subrepresentation of $\etcoh{i}{\X_{\bar{\Fi}_p}}$, it will inherit these properties. Thus we can use the arguments of \cite{EB} to conclude from equation (\ref{eq-poly}) that the semisimplification of $\etcohl{i}{\X_{\bar{\Q}}}$ has dimension zero if $i$ is odd and is isomorphic to $\Q_l(-i/2)^{P_{i/2}}$ if $i$ is even. Then the claim follows, analogously as in \cite{EB}, from the fact that $\etcohl i {\X_{\bar{\Q}}}$ is potentially semistable. This is, again, a consequence of the fact that it is a subrepresentation of $\etcoh i {\X_{\bar\Q}}$, which is shown to be potentially semistable in \cite{EB}.
\qed

Note that the moduli spaces $\Mmb gn$ are indeed proper and smooth Deligne-Mumford stacks defined over $\Z$, of pure relative dimension and with an action of $\s_n$ (see \cite{Deligne-Mumford} and \cite[Theorem 2.7]{Knudsen}).

\subsection{The mixed Hodge structures} \label{count-to-betti}
Using comparison theorems we get a theorem corresponding to Theorem~\ref{thm-count-adic}, but for Betti cohomology. This new theorem is an equivariant version of Corollary 5.3 of \cite{EB}.

Let $\X$ be an algebraic stack (or simply a scheme) defined over $\C$ together with an action of $\s_n$. Analogously to the case of $l$-adic cohomology in Section \ref{etGK}, we denote by $H_{c,\lambda}^i(\X;\Q)$ the sum of all copies of the irreducible representation of $\s_n$ indexed by the partition $\lambda$ that appear in the Betti cohomology groups with compact support $H_c^i(\X;\Q)$. We denote by $\Eul_\lambda(\X)$ the Euler characteristic of $H_{c,\lambda}^\pu(\X;\Q)$ in $\mathrm K_0(\mathsf{MHS}_\Q)$.

Since all the isomorphisms used in the proof of Corollary 5.1 of \cite{EB} respect the action of the symmetric group we have a proof of the following theorem.

\begin{theorem} \label{thm-count-hodge}
In the situation of Theorem \ref{thm-count-adic} suppose furthermore that the coarse moduli space of the stack $\X_{\Q}$ is the quotient of a smooth projective $\Q$-scheme by a finite group. 

Then for all partitions $\lambda$ of $n$ and for all $i \geq 0$, there is an isomorphism of pure $\Q$-Hodge structures
$$H^i_{c,\lambda}( \X(\C) , \Q) \simeq \begin{cases} 0 &
  \text{if i is odd} \\ \Q(-i/2)^{b_{i/2}} & \text{if i is even}
  \end{cases}$$
where the left hand side is equipped with the canonical Hodge structure of \cite{Deligne}.
\end{theorem}

Note that the moduli spaces $\Mmb gn$ fulfill the additional condition of Theorem \ref{thm-count-hodge} (see \cite{Boggi-Pikaart}).

\subsection{The main theorem} \label{count-gives}
The following theorem summarizes what cohomological information we can get from polynomial counts of the points of the spaces $\Mm gn$.

\begin{theorem} \label{thm-count-stable}
Assume that, for all $\Mm gn$ with $(g,n) \in \mathcal{D}_{G,N}$ and for all partitions $\lambda$ of $n$, equation (\ref{cond-adic}) of Theorem \ref{thm-count-adic} is fulfilled for some set $S$ and polynomial $P_{\lambda,g,n}(t)$.

Then the following holds for all $(g,n) \in \mathcal{D}_{G,N}$ and all partitions $\lambda$ of $n$.

For $\Mmb gn$ equation (\ref{cond-adic}) of Theorem \ref{thm-count-adic} is fulfilled for some set $S$ and polynomials $Q_\lambda(t)$. Hence, both Theorem~\ref{thm-count-adic} and Theorem~\ref{thm-count-hodge} hold for $\Mmb gn$.

Moreover, 
\begin{itemize}
\item[(i)] $\etEull(\Mm gn \otimes_{\Z} {\bar{\Q}}) =P_{\lambda,g,n}([\Q_l(-1)]) \in  \mathrm{K}_0(\mathsf{Gal}_{\Q_l})$
\item[(ii)] $\Eull(\Mm gn (\C))=P_{\lambda,g,n}(\Ll) \in \mathrm K_0(\mathsf{MHS}_\Q) $.
\end{itemize}
\end{theorem}

\begin{proof}
This is proved by induction on $(g,n)$. The starting case is $\Mm 03$ for which the result is clear since we can apply Theorems \ref{thm-count-adic} and \ref{thm-count-hodge} to $\Mm 03 = \Mmb 03$. 

Say that the theorem holds for $\mathcal{D}_{g,n} \setminus \{(g,n)\}$. In the stratification of $\Mmb gn$, $\Mm gn$ is the open part and it does not contribute to the cohomology of the boundary $\partial\Mm gn=\Mmb gn\setminus\Mm gn$. From the additivity of Euler characteristics we get an equality
\begin{equation} \label{eq-ind}
\etEull(\Mmb gn \otimes_{\Z} \bar{\Q}) = \etEull(\Mm gn \otimes_{\Z} \bar{\Q}) + \etEull(\partial \Mm gn \otimes_{\Z} \bar{\Q}). 
\end{equation}
Thus, using Theorem \ref{thmGK} and Section \ref{etGK} we can express $\etEull(\partial\Mm gn \otimes_{\Z} \bar{\Q})$ in terms of $\etEulu(\Mm{\hat{g}}{\hat{n}} \otimes_{\Z} \bar{\Q})$ for partitions $\mu$ of $n$ and $(\hat{g},\hat{n}) \in \mathcal{D}_{g,n} \setminus \{(g,n)\}$. Since, by induction, $\etEulu(\Mm{\hat{g}}{\hat{n}} \otimes_{\Z} {\bar{\Q}})$ is known for these indices and is a sum of elements $[\Q_l(-j)]$, we can determine $\etEull(\partial\Mm gn \otimes_{\Z} \bar{\Q})$ and it will also be a sum of elements $[\Q_l(-j)]$.

Equation (\ref{eq-ind}) takes place in the category $\mathrm{K}_0(\mathsf{Gal}_{\Q_l})$ and thus we can take the trace of Frobenius on both sides. By the assumption together with Section \ref{count-to-trace}, for almost all finite fields $k$ we can compute $\mathrm{Tr}(F_q,\etEull(\Mm gn \otimes_{\Z} {\bar{k}}))$, which is equal to $\mathrm{Tr}(F_q,\etEull(\Mm gn \otimes_{\Z} {\bar{\Q}}))$, see \cite{EB}. Just as $\mathrm{Tr}(F_q,\etEull(\partial\Mm gn\otimes_{\Z} \bar{\Q}))$, it will be polynomial in $\abs{k}=q$.

Hence we can use equation (\ref{eq-ind}) to compute $\mathrm{Tr}(F_q,\etEull(\Mmb gn \otimes_{\Z} {\bar{\Q}}))$ and we know that the answer will be polynomial in $q$. Theorems \ref{thm-count-adic} and \ref{thm-count-hodge} are therefore applicable for $\Mmb gn$.

Thus we know both $\etEull(\Mmb gn \otimes_{\Z} {\bar{\Q}})$ and $\etEull(\partial\Mm gn\otimes_{\Z} \bar{\Q})$, which together with equation (\ref{eq-ind}) gives (i). But equation (\ref{eq-ind}) also holds in the category of mixed Hodge structures and therefore we can use the same argument to conclude (ii). Here the induction ends.
\end{proof}

\begin{remark} \label{rem-genus01}
In \cite{Lehrer}, Kisin and Lehrer made an equivariant count of the number of points of $\Mm 0N$ for $N \geq 3$ and these numbers all fulfill the hypotheses of Theorem~\ref{thm-count-stable}.

The first author has made an equivariant count of the numbers of points of $\Mm 1N$ for $N\leq 6$, which also were found to be polynomial. It was achieved by applying a method similar to that described in Section \ref{quartics}, and based on the embedding of a genus one curve with a marked point $P$, via the divisor $3P$, in $\Pb^2$. This computation agrees with Getzler's results on the Hodge Euler characteristic of $\Mm 1N$ for $N \geq 1$ (see \cite{G-res}).
\end{remark} 

\section{Counting points over finite fields}\label{PartIII}
For all finite fields $k$, with the possible exception of a finite number of characteristics, we wish to make an  $\s_n$-equivariant count (for the definition see Section \ref{count-to-trace}) of the number of points defined over $k$ of $\Mm{g}{n} \otimes_{\Z} \bar{k}$ for all $(g,n) \in \mathcal{D}_{4,0} \setminus \{(4,0)\}$. 

Following Section \ref{count-to-trace}, we will present results of equivariant counts of the number of points in the form of traces of Frobenius on the $\s_n$-equivariant Euler characteristics in the category of Galois representations
$$\etEuls{n}(\Mm gn \otimes_{\Z} \bar{k}) := \sum_{\lambda \, \vdash n} \chi_{\lambda}(id)^{-1} \etEull\bigl(\Mm gn\otimes_{\Z} \bar{k} \bigr) s_{\lambda} \in \mathrm K_0(\mathsf{Gal}) \otimes \Lambda_n.$$

\subsection{Preparation}\label{what-to-count}
In the following we will always assume that $\Mm{g}{n}$ is defined over $\bar k$, and we denote by $M_{g,n}$ its coarse moduli space. Recall that $F_q$ denotes the geometric Frobenius belonging to $k$.

\begin{definition}
An $n$-tuple of distinct points $(p_1, \ldots, p_n)$ is called a \emph{conjugate $n$-tuple} if $F(p_i)=p_{i+1}$ for $1 \leq i \leq n-1$ and $F(p_n)=p_{1}$.

If $\lambda = (\lambda_1, \ldots, \lambda_\nu)$ then a $\abs{\lambda}$-tuple of distinct points $(p_1, \ldots, p_{\abs{\lambda}})$ is called a \emph{$\lambda$-tuple} if $(p_{\sum_{i=1}^{j-1} \lambda_i+1},p_{\sum_{i=1}^{j-1}\lambda_i+2}, \ldots, p_{\sum_{i=1}^{j} \lambda_i})$ is a conjugate $\lambda_j$-tuple for every $1 \leq j \leq \nu$. \end{definition}

From the definition of the $\s_n$-action on $\Mm gn$ we see that if $\sigma \in \s_n$ has cycle type $\lambda$ then the number of fixed points of $F \cdot \sigma$ acting on $M_{g,n}$ is equal to the number of $\bar k$-isomorphism classes of curves $C$ in the category $\M{g}(k)$ together with a $\lambda$-tuple of points $(p_1, \ldots, p_n)$ lying on the curve. Here we use the fact that a $\bar{k}$-isomorphism class of $\M{g}(\bar{k})$ is fixed by Frobenius if and only if it comes from an element of $\M{g}(k)$ (see \cite[Lemma 10.7.5]{Katz-Sarnak}).

The Lefschetz fixed point theorem (see for instance \cite[Theorem V.2.5]{Milne}) shows that for all $C$ in $\M{g}(k)$ and $m \geq 1$ 
$$ \abs{C(k_m)} = \sum_i (-1)^i \mathrm{Tr}(F_q^m,\etcoh{i}{C_{\bar{k}}}) = 1+ q^m-a_m(C) $$ 
where $a_{m}(C):=\mathrm{Tr} ( F_q^m,\etcoh{1}{C_{\bar{k}}})$. 
Note that the number of conjugate $m$-tuples on a curve $C$ equals $$\sum_{d | m} \mu(m/d) \abs{C(k_m)}$$ where $\mu$ is the M\"obius function. Instead of counting the number of $\bar k$-isomorphism classes of the curves described above, one can count $k$-isomorphism classes but then with weight the reciprocal of the number of $k$-automorphisms of the curve (see \cite[Proposition 5.1]{Geer} or \cite[Lemma 10.7.5]{Katz-Sarnak}).
Together this shows that 
$$\abs{M_{g,n}^{\sigma \cdot F_q}} = \sum_{C \in [\M{g}(k)]} \frac{Q_{\lambda}(q,a_1(C), \ldots , a_n(C))}{\abs{\Aut_k(C)}} $$ 
where $Q_{\lambda}(x_0,x_1, \ldots, x_n)$ is a polynomial with coefficients in $\Z$, and the sum runs over all isomorphism classes of the category $\M g(k)$. If we assign degree $i$ to the indeterminate $x_i$, then $Q_{\lambda}(x_0,x_1, \ldots, x_n)$ has a single monomial of the highest degree $\abs{\lambda}$, namely $(-1)^{\ell(\lambda)}x_{\lambda_1} \cdots x_{\lambda_{\nu}}$. 

Thus, if $\abs{M_{g,m}^{\tau \cdot F_q}}$ is known for all $m < n$ and for all $\tau \in \s_m$, the only part missing to compute $\abs{M_{g,n}^{\sigma \cdot F_q}}$ is $$\sum_{C \in [\M{g}(k)]} \frac{1}{\abs{\Aut_k(C)}} \cdot \prod_{i=1}^\nu a_{\lambda_i}(C).$$ 

\subsection{Hyperelliptic curves}\label{counthyp}
In this section we assume that $k$ is a finite field of odd characteristic. 
Let $\Hhm{g}{n}$ denote the subset of $\Mm{g}{n}$ of hyperelliptic curves. We wish to make an equivariant count of the number of points defined over $k$ of the space $\Hh{g,n}$ for all $g \in\{2,3\}$ and $n \leq 8-2g$.  

In view of the results in section \ref{what-to-count}, it is equivalent to compute $$\sum_{C \in [\Hh{g}(k)]} \frac{1}{\abs{\Aut_k(C)}} \cdot \prod_{i=1}^\nu a_{\lambda_i}(C)$$ for all partitions $\lambda=(\lambda_1, \ldots, \lambda_{\nu})$ of weight $n$.

Viewing the hyperelliptic curves of genus $g \geq 2$ defined over $k$ as double covers of $\Pb^1$, we can write them in the form $y^2=f(x)$, where $x$ is a local coordinate on $\Pb^1$ and $f$ is in $P_g$, the set of square-free polynomials of degree $2g+1$ or $2g+2$ with coefficients in $k$. It follows that we can sum over the set of curves corresponding to the elements of $P_g$ and then divide by the number of elements of the group of isomorphisms defined over $k$ of these curves, that is, by $\abs{\mathrm{GL}_2(k)}$.

If we denote by $C_f$ the curve corresponding to the polynomial $f \in P_g$, then by the Lefschetz fixed point theorem  we have
$$a_m(C_f) = -\sum_{\alpha \in   \Pb^1(k_m)} \chi_{2,m} \bigl( f(\alpha) \bigr),$$ 
where $\chi_{2,m}$ is the quadratic character of the field $k_m$. Recall that the quadratic character of a field is the function mapping 0 to 0, all non-zero squares of the field to 1, and all other elements to $-1$. We conclude that
$$\sum_{C \in [\Hh{g}(k)]} \frac{1}{\abs{\Aut_k(C)}} \cdot \prod_{i=1}^\nu a_{\lambda_i}(C) = \frac{(-1)^{\abs{\lambda}}}{\abs{\mathrm{GL}_2(k)}}  \sum_{f \in P_g} \prod_{i=1}^{\nu} \Bigl(\sum_{\alpha \in   \Pb^1(k_{\lambda_i})} \chi_{2,\lambda_i} \bigl( f(\alpha) \bigr) \Bigr).
$$
This sum splits into sums of the form
 $$u_{g}(\alpha_1,\dots,\alpha_\nu):= \sum_{f \in   P_g} \prod_{i=1}^\nu \chi_{2,\lambda_i} \bigl( f(\alpha_{i}) \bigr),$$ 
for every choice of points $\alpha_i \in \Pb^1(k_{\lambda_i})$. Note that if $n$ is odd, then all these sums are equal to zero.

\begin{example}
We consider the case $\lambda=(1,1)$ and $(\alpha_1,\alpha_2)=(0,\infty)$. We want to compute $u_g:=u_g(0,\infty)$, for all $g\geq0$. If $f$ is a monic polynomial, representing a hyperelliptic curve, then $f(\infty)=0$ if it is of odd degree and $f(\infty)=1$ if it is of even degree. Thus $u_g$ is equal to the sum of $\chi_{2,1} \bigl( f(\alpha_1) \bigr)$ multiplied with $q-1$, over all monic square-free polynomials $f$ of degree $2g+2$. Since there are as many non-zero squares as non-squares in $k$, the sum of $\chi_{2,1} \bigl( h(\alpha_1) \bigr)$ over all monic polynomials $h$ is zero.
Hence we get the following recursion formula  for $g \geq 0$, 
\begin{multline*} 
0 = u_g + u_{g-1} (q-1) + u_{g-2} (q-1) q + \ldots + u_{0} (q-1) q^{g-1} + (q-1)^2 q^{g}
\end{multline*} from which we conclude that $u_g(0,\infty)=-(q-1)^2$ for all $g\geq 0$. 
\end{example}

The example above can be generalized in the following sense. Fix a $\lambda$-tuple $(\alpha_1,\dots,\alpha_\nu)$ of points, which we for simplicity assume to be in $\Ab^1$, and set $u_g:=u_g(\alpha_1,\dots,\alpha_\nu)$. Then we can determine $u_g$ recursively, once we know the integers 
$$U_{g}:= \sum_{h}\prod_{i=1}^\nu \chi_{2,\lambda_i} \bigl( h(\alpha_{i}) \bigr),$$ 
where the sum is over all polynomials $h$ of degree $2g+1$ or $2g+2$.

Define $b_j$ to be the number of monic polynomials $l$ of degree $j$ such that $l(\alpha_i)$ is non-zero for all $i$. Any polynomial $h$ has a unique decomposition of the form $h = f \cdot l^2$ where $f$ is a square-free polynomial and $l$ is a monic polynomial. Thus for any point $\alpha \in \Ab^1(k_m)$, $\chi_{2,m} \bigl( h(\alpha) \bigr) = \chi_{2,m} \bigl(f(\alpha) \bigr)$ as long as $l(\alpha)$ is non-zero. From this we can conclude that 
\begin{equation}\label{induc}
U_g = (q-1)b_{g+1} + \sum_{i=0}^{g} b_{i} u_{g-i} .
\end{equation}

If $h$ is a polynomial of degree greater than $n-1$, there is a one to one correspondence between the set of possible values of the first $n$ coefficients of $h$ and the set of possible values of $h(\alpha_i)$ for $1 \leq i \leq \nu$. 
Using this, and the fact that half of the non-zero elements of any finite field of odd characteristic are squares and half are non-squares, we can determine $U_g$ for all $g$ satisfying $2g+1\geq n$.

This shows that, if we know the values of $u_g$ for all $g$ such that $2g+1 < n$, we can use \eqref{induc} to determine their value for all possible $g$.

In \cite{Jonas1} the method described above, which follows a suggestion by Nicholas M. Katz, together with the fact that if $C$ has genus zero then $a_m(C)=0$ for all $m$, is used to make an equivariant count of the number of points defined over $k$ of $\Hhm{g}{n}$ for any value of $g$ and all $n \leq 7$. The results, for the cases relevant here, are given in Table~\ref{tabH}. In particular, in view of Theorem~\ref{thm-count-stable}, the equivariant count of points of $\Hhm 24$ implies Theorem~\ref{char24}.

\def\quod{\hskip 0.5em\relax }
\begin{table}\caption{\label{tabH} Equivariant counts of points of some $\Hhm gn$.}
\vbox{\bigskip}
\centerline{
\vbox{
\offinterlineskip
\hrule
\halign{&\vrule#& \quod #\hfil \strut \quod \cr
height2pt&\omit&&\omit& \cr
&$n$&&$\mathrm{Tr}\bigl(F_q,\etEuls{n}(\Hhm{2}{n} \otimes_{\Z} \bar{k})\bigr)$&\cr
height2pt&\omit&&\omit&\cr
\noalign{\hrule}
height2pt&\omit&&\omit&\cr
&$0$&&$q^3  s_0$&\cr
height2pt&\omit&&\omit&\cr
&$1$&&$(q^4+q^3)  s_1$&\cr
height2pt&\omit&&\omit&\cr
&$2$&&$(q^5+q^4-1)  s_2 + q^4  s_{1^2}$&\cr
height2pt&\omit&&\omit&\cr
&$3$&&$(q^6+q^5-q^4-q)  s_3 + (q^5-q)  s_{21}$&\cr
height2pt&\omit&&\omit&\cr
&$4$&&$(q^7+q^6-q^5-q^4-q^2+q+1)  s_4 + (q^6-q^5-q^2)  s_{31}$&\cr
height2pt&\omit&&\omit&\cr
&\omit&&$- q^2  s_{2^2} + (q^2-1)  s_{21^2}$&\cr
height2pt&\omit&&\omit&\cr
\noalign{\hrule}
height2pt&\omit&&\omit&\cr
&$n$&&$\mathrm{Tr}\bigl(F_q,\etEuls{n}(\Hhm{3}{n}\otimes_{\Z} \bar{k})\bigr)$&\cr
height2pt&\omit&&\omit&\cr
\noalign{\hrule}
height2pt&\omit&&\omit&\cr
&$0$&&$q^5  s_0$&\cr
height2pt&\omit&&\omit&\cr
&$1$&&$(q^6+q^5)  s_1$&\cr
height2pt&\omit&&\omit&\cr
&$2$&&$(q^7+q^6)  s_2 + (q^6-1)   s_{1^2}$&\cr
height2pt&\omit&&\omit&\cr
} \hrule}}
\end{table}

\subsection{Quartic curves} \label{quartics}
Let $\mathcal{Q}_n$ be the complement of $\Hhm{3}{n}$ in $\Mm{3}{n}$. Using the canonical embedding we can identify $\mathcal{Q}_n$ with the moduli space of $n$-pointed plane non-singular quartic curves. By the results in Section \ref{what-to-count} we find that an equivariant count of the number of points of $\mathcal{Q}_n$ is the same as an equivariant count of the number of plane $n$-pointed non-singular quartic curves defined over $k$ divided by the number of elements, defined over $k$, of the group of isomorphisms acting on these curves, which is equal to $\mathrm{PGL}_3(k)$.

From now on all curves mentioned in this section will be assumed to be plane curves.

Fix a partition $\lambda$ of weight $n$. We want to compute the fixed points of $\sigma \cdot F_q$ on the coarse moduli space of $\mathcal{Q}_n$, where $\sigma$ is any permutation with cycle type $\lambda$. By Section~\ref{what-to-count} we want, in other words, to compute the sum, over all $\lambda$-tuples $P$ of points in the plane, of the number of non-singular quartic curves that contain $P$.

\begin{definition}
Let us identify the space of quartic curves defined over $k$ with $\Pb^{14}(k)$. 
For every $\lambda$-tuple $P$ of points, we denote by $L_P$ the linear subspace of $\Pb^{14}(k)$ of curves that contain $P$. 

A set $S$ of $m$ distinct points in the plane is called an \emph{unordered $\mu$-tuple} if there is a $\mu$-tuple $(p_1,\dots,p_m)$ such that $S=\{p_1,\dots,p_m\}$. For any unordered $\mu$-tuple $S$ of points in the plane, we define $L_{P,S}$ to be the linear subspace of $L_P$ of curves that have singularities at the points of $S$.
\end{definition}

The locus of singular curves in $L_P$ is the union of all linear spaces $L_{P,S}$ for every $m\geq 1$ and unordered $m$-tuple $S$. Say that we want to use the sieve principle to count the number of elements of this union. Then we would sum the numbers $(-1)^{i+1}  \abs{L_{P,S_1} \cap \ldots \cap L_{P,S_i}}$, for each $i \geq 1$ and for each unordered choice of distinct sets $S_1, \ldots, S_i$ where $S_j$ is an unordered $m_j$-tuple. If this procedure terminated, every singular curve in $L_P$ would have been counted exactly once. Hence, taking $\abs{L_P}$ minus the resulting number would give the number of non-singular curves in $L_P$. Since there are curves with infinitely many singularities, namely the non-reduced ones, this procedure will not stop. Instead we choose to stop it after making sure that all curves with up to two singularities have been removed precisely once. That is, we will compute,
\begin{multline}
s_{2,P}:=\abs{L_P}-\sum_{\alpha \in \Pb^2(k)} \abs{L_{P,\{\alpha\}}} + \sum_{\alpha_1 \neq \alpha_2 \in \Pb^2(k)} \frac{1}{2}\abs{L_{P,\{\alpha_1,\alpha_2\}}} \\- \sum_{\beta \in \Pb^2(k_2) \setminus \Pb^2(k)} \frac{1}{2}\abs{L_{P,\{\beta,F(\beta)\}}}.
\end{multline}

\begin{example} 
Suppose that $P$ consists of the point $p_1$. In this case $L_P$ has dimension $13$. If $p_1$ is not contained in $S$, the subspace $L_{P,S}$ will have the expected dimension; if this is not the case, the dimension of $L_{P,S}$ will be one more than expected. Hence, we have to distinguish the following five types of singular curves passing through $P$: curves singular in $p_1$, singular in a point different from $p_1$, singular in a conjugate $2$-tuple, singular in $p_1$ and another $k$-point and finally those singular in two $k$-points both different from $p_1$. Note that these cases are not mutually exclusive. By choosing the appropriate signs we get, 
\begin{multline*}
s_{2,P} = \abs{\Pb^{13}(k)} - \abs{\Pb^{11}(k)} - (q^2+q)
 \abs{\Pb^{10}(k)} - \frac{1}{2}(q^4-q)  \abs{\Pb^{7}(k)} +
\\ + (q^2+q)  \abs{\Pb^{8}(k)} + \frac{1}{2}(q^2+q) 
(q^2+q-1)  \abs{\Pb^{7}(k)}. 
\end{multline*}
\end{example}

After having computed the sum $\sum_{P} s_{2,P}$ where $P$ runs over all $\lambda$-tuples, we need to amend for the curves with more than two singularities to obtain the sum, over all $\lambda$-tuples $P$, of the number of non-singular curves in $L_P$. 

The choice of stopping the sieve procedure after two singularities was made after weighing the difficulty of computing the dimensions of the linear subspaces $L_{P,S}$ against the difficulty of finding the more specific information we need on the singular curves that have not been removed precisely once in the sieve procedure.

The information we need on the curves that have not been removed precisely once is the following. For any partition $\mu$ of weight greater than two, let $t_{\mu}$ be the sum over all choices of $\lambda$-tuples $P$ and unordered $\mu$-tuples $S$ of the number of curves that contain $P$ and that have singularities at the points of $S$, and nowhere else. Each number $t_{\mu}$ will then have to be removed or added to $\sum_{P} s_{2,P}$ a suitable number of times, so that the curves of this kind also will have been removed precisely once.

All quartic curves with at least three singularities are such that all components (or their normalizations) have genus zero, with the only exception of those consisting of a non-singular cubic together with a transversal line. To find all the numbers $t_{\mu}$ we distinguish types of these quartic curves defined over $k$ according to the following information: 
\begin{itemize}
\item degrees and multiplicities of the irreducible components; 
\item over which fields the components are defined; 
\item the number of singularities of each irreducible component;
\item over which fields the singularities are defined;
\item the number of intersection points of any set of irreducible components;
\item over which fields the intersection points are
defined. \end{itemize}

Hence if we know the type of a curve we only need to find the number of points of its irreducible components to deduce its contribution to $t_{\mu}$.

Since irreducible curves of degree one or two are isomorphic to the projective line their number of points is easy to compute.

\begin{example}
Let us consider quartic curves defined over $k$ consisting of a conic together with one tangent line and one transversal line intersecting the conic in two points defined over $k$ distinct from the tangency point of the other line. In this case all components and intersection points will be defined over $k$. Furthermore all components are rational curves, so the number of points defined over $k_i$ of such a curve equals $3  \abs{\Pb^1(k_i)}-4$. The number of such curves is easily computed to be $\abs{\mathrm{PGL}_3(k)}/2$.
\end{example}

In case the normalization of a singular component has genus zero we only need to find how many points there are in the inverse image of the singularities under the normalization map and over which fields they are defined, to be able to compute its number of points.

\begin{example}
There are $(q-2)\abs{\mathrm{PGL}_3(k)}/6$ quartic curves defined over $k$ consisting of a cuspidal cubic together with a transversal line intersecting the cubic in three points defined over $k$. The fact that the singularity of the cubic is a cusp means that the inverse image of the singularity under the normalization map is only one point which is defined over the same field as the singularity. This shows that curves of this type have $2\abs{\Pb^1(k_i)}-3$ points defined over $k_i$.
\end{example}

Non-singular cubic curves have genus one, and hence we cannot use the same approach as above. But pointed non-singular cubic curves can be counted in the same fashion as we are counting the pointed non-singular quartic curves, see also Remark \ref{rem-genus01}. Hence after adding the transversal line we can also find the contribution of these curves to $t_{\mu}$.

This method was used in \cite{Jonas2} and gave the equivariant answers for the number of points defined over $k$ of $\mathcal{Q}_{n}$ for all $n$ less than or equal to six. The results for the cases relevant here can be found in Table~\ref{tabM3}. In particular, the results for $\mathcal Q_{2}$, together with the results on $\Hhm 32$ given in Table~\ref{tabH}, imply Theorem~\ref{char32}.

\begin{table}\caption{\label{tabM3}Equivariant counts of points of $\mathcal Q_{n}$ for $n\leq2$.}
\vbox{\bigskip}
\centerline{
\vbox{
\offinterlineskip
\hrule
\halign{&\vrule#& \quod #\hfil \strut \quod \cr
height2pt&\omit&&\omit& \cr
&$n$&&$\mathrm{Tr}\bigl(F_q,\etEuls{n}(\mathcal{Q}_n\otimes_k \bar{k})\bigr)$&\cr
height2pt&\omit&&\omit&\cr
\noalign{\hrule}
height2pt&\omit&&\omit&\cr
&$0$&&$(q^6+1)  s_0$&\cr
height2pt&\omit&&\omit&\cr
&$1$&&$(q^7+q^6+q+1)  s_1$&\cr
height2pt&\omit&&\omit&\cr
&$2$&&$(q^8+q^7-q^5+q^2+q)  s_2 + (q^7+q+1)  s_{1^2}$&\cr
height2pt&\omit&&\omit&\cr
} \hrule}}
\end{table}

\section{Rational cohomology of geometric quotients}\label{PartII}

If $g$ and $n$ are small enough, it is possible to determine exactly what the rational cohomology of $\Mm gn$ is. This applies at least for the moduli space of curves of genus 3 and 4. The reason is that in these cases the canonical map provides a very simple description of each such curve. This allows to divide the coarse moduli space into locally closed subschemes which are geometric quotients of the complement of a discriminant in a complex vector space, by the action of a reductive group. In this section, we explain how this description applies to moduli spaces of ($n$-pointed) non-singular curves of genus 3 or 4. Moreover, following \cite{PS}, we explain which relation holds between the cohomology of the complement of the discriminant and the cohomology of the geometric quotient. 

Note that, throughout this section, we will always consider not the moduli stacks $\Mm gn$, but the underlying coarse moduli spaces. Since the rational cohomology of the two coincides, we will abuse notation and not distinguish between the two concepts.

Throughout this section and the following one, all results on the cohomology (or the Borel-Moore homology) will be given by means of \emph{Poincar\'e-Serre polynomials} (following \cite{Looij}). If $H^\pu$ is a rational mixed Hodge structure, the Poincar\'e-Serre polynomial of $H^\pu$ is defined as the polynomial in $\Z[t,u,u^{-1}]$ such that the coefficient of $t^iu^{j}$ is the dimension of the weight $j$ subquotient of $H^i$. The Poincar\'e-Serre polynomial of a complex variety $Z$ is the Poincar\'e-Serre polynomial of its rational cohomology. Note that in all cases we will consider, the mixed Hodge structures will be sums of (rational) Tate Hodge structures. 
In this special case, giving the Poincar\'e-Serre polynomial is equivalent to giving the whole rational mixed Hodge structure.

\subsection{Moduli spaces of smooth subvarieties}\label{introLHic}

Let $Z\subset\Pp N$ be a complex projective variety, and $L$ a vector bundle on $Z$. We define $V$ to be a linear subspace of the space of sections of $L$. Inside $V$ there is a closed subvariety $\Sigma$ of sections whose zero scheme is not smooth of the expected dimension. We call $\Sigma$ the \emph{discriminant}. In the following, we will always assume that $\Sigma$ is of pure complex codimension 1 in $V$.

Suppose further that there is a complex affine algebraic group $G$ acting on $Z$ and $L$, and that this action induces an action on $V$.

We are interested in studying the quotient of $X:=V\setminus\Sigma$ by the action of $G$ (if such a quotient exists), and specifically in determining its rational cohomology. In many examples the geometric quotient $X/G$ can be interpreted as the moduli space of the smooth varieties defined by the vanishing of elements in $V$.

Moreover, for every $n\geq0$ we can consider the incidence correspondence
$$\mathcal I_n:= \{(v,p_1,\dots,p_n)\in X\times\F nZ: v(p_i)=0\ \forall i=1,\dots,n\},$$
where $\F nZ$ denotes the space of ordered configurations of $n$ distinct points on $Z$.
The action of $G$ on $Z$ clearly induces an action on $\mathcal I_n$ and if  the geometric quotient $X/G$ exists, so does $\mathcal I_n/G$. 

We recall Peters-Steenbrink's generalization of the Leray-Hirsch theorem:

\begin{theorem}[\cite{PS}]\label{LeHi}
Let $\varphi: X\rightarrow Y$ be a geometric quotient for the action of a connected group $G$, such that for all $x\in X$ the connected component of the stabilizer of $x$ is contractible. Consider the orbit inclusion
$$\begin{matrix}\rho:&G&\longrightarrow&X\\&g&\longmapsto&gx_0,\end{matrix}$$
where $x_0\in X$ is a fixed point.

Suppose that for all $k>0$ there exist classes $e_1^{(k)},\dots,e_{n(k)}^{(k)}\in\coh kX$ that restrict to a basis for $\coh kG$ under the map induced by $\rho$ on cohomology. Then the map $$a\otimes\rho^*(e_i^{(k)})\longmapsto\varphi^*a\cup e_i^{(k)}$$
extends linearly to an isomorphism of graded linear spaces
$$\coh\pu Y\otimes\coh\pu G \cong \coh\pu X$$
that respects the rational mixed Hodge structures of the cohomology groups.
\end{theorem}

\begin{remark}\label{casiLH}
The hypotheses of Theorem~\ref{LeHi} are apparently very natural in all cases where $G$ is a reductive group and $X=V\setminus\Sigma$ is the complement of the discriminant.

 In particular, they are known to hold in the following cases:\begin{enumerate}
\item\label{hs}\emph{Moduli spaces of smooth hypersurfaces} \cite[Theorem 1]{PS}\\
$Z=\Pp N$, $V$ is the vector space of sections of $L=\mathcal{O}_{\Pp N}(d)$ ($d\geq 3$) and $G=\GL(N,\C)$.
\item\label{4co}\emph{Moduli space of smooth curves on a quadric cone} \cite[4.1]{OT}\\
$Z=\Pp{}(1,1,2)$, $V$ is the vector space of sections of $L=\mathcal{O}_Z(6)$ and $G$ is the automorphism group of the graded ring $\C[x,y,z]$ where $\deg x=\deg y=1,\ \deg z=2$. This can also be generalized to curves of even degree $\geq 4$ on a quadric cone.
\item\label{4ge}\emph{Moduli space of smooth $(3,3)$-curves on $\Pp1\times\Pp1$} \cite[3.1]{OT}\\
$Z=\Pp1\times\Pp1$, $V$ is the vector space of sections of $L=\mathcal{O}_{\Pp1}(3)\otimes\mathcal{O}_{\Pp1}(3)$ and $G$ is the connected component of the identity of the automorphism group of the $\Z\times\Z$-graded ring $\C[x_0,x_1;y_0,y_1]$. Note that this result can be easily generalized to smooth $(m,n)$-curves on $\Pp1\times\Pp1$, with $m,n\geq2$.
\end{enumerate}
\end{remark}

\begin{theorem}\label{inco}
Suppose $X$ satisfies the hypotheses of Theorem~\ref{LeHi} for the action of an affine algebraic group $G$. Then the action of $G$ on $\mathcal I_n$ satisfies the hypotheses of Theorem~\ref{LeHi}, for every $n\geq0$.
\end{theorem}

\proof
Consider the natural projection $\pi_n: \mathcal I_n\rightarrow X$. The map $\pi_n$ is equivariant with respect to the action of $G$. Then the claim follows from the fact that, for an appropriate choice of the base points, the orbit inclusion of $G$ in $X$ is the composition of $\pi_n$ with the orbit inclusion of $G$ in $\mathcal I_n$.\qed

\subsection{Moduli spaces of curves of genus three and four}\label{M3}

The considerations above can be applied directly to the study of the rational cohomology of moduli spaces of curves of genus 3 and 4, with or without marked points. Indeed, every non-hyperelliptic curve of genus three has a plane quartic curve as its canonical model (\cite[IV.5.2.1]{Hartshorne}). This means that we are in the situation of Remark~\ref{casiLH}, case \ref{hs}, with $N=2$ and $d=4$. In the notation of that case, we have that
$$\coh\pu{\mathcal Q_{n}}\otimes\coh\pu{\GL(3)}\cong \coh\pu{\mathcal I_n},$$
where $\mathcal Q_{n}$ is the complement in $\Mm3n$ of the hyperelliptic locus $\Hh{3,n}$.
Hence, determining the cohomology of $\mathcal Q_{n}$ is equivalent to determining that of $\mathcal I_n$. 

We are interested in determining the cohomology of $\mathcal I_n$ for the first values of $n$. For $n=0$, we have $\mathcal I_0=X=V\setminus\Sigma$. The cohomology of this space has been computed by Vassiliev in \cite{Vart}. Vassiliev's results imply that the Poincar\'e-Serre polynomial of $\mathcal I_0$ is equal to
$$(1+\ut2{})(1+\ut43)(1+\ut65)(1+\ut{12}6).$$

By Theorem~\ref{inco}, this implies that the Poincar\'e-Serre polynomial of $\mathcal Q_{0}$ is $1+\ut{12}6$. Since the hyperelliptic locus has the rational cohomology of a point, this yields an alternative proof of the following result of Looijenga.

\begin{theorem}[{\cite[4.7]{Looij}}]
The Poincar\'e-Serre polynomial of $\M3$ is $1+\ut22+\ut{12}6$.
\end{theorem}

Consider next the case $n=1$. The incidence correspondence $\mathcal I_1$ has a natural forgetful morphism $\pi_1: \mathcal I_1\longrightarrow \Pp2$. The map $\pi_1$ is a locally trivial fibration, and the fibre $F_1$ is isomorphic to the space of non-singular homogeneous polynomials of degree 4, vanishing at a fixed point $p\in\Pp2$. The rational cohomology of this space is computed with Vassiliev-Gorinov's method in Section~\ref{VG} (see \eqref{I1} at page \pageref{I1}). In particular, its Poincar\'e-Serre polynomial is
$$(1+\ut2{})(1+2\ut43+\ut86+\ut{12}6+2\ut{16}9+\ut{20}{12}).$$

Note the the rational cohomology of $F_1$ is the tensor product of the cohomology of $\C^*$ and the cohomology of the projectivization of $F_1$. The Leray spectral sequence associated to $\pi_1$ is given in Table~\ref{tabM31}, where for the sake of simplicity we have already divided out the cohomology of $\C^*$.

\begin{table}[t]\caption{Leray spectral sequence converging to $\coh\pu{\mathcal Q_1}\otimes\coh\pu{\PGL(3)}$. \label{tabM31}}
$$\begin{array}[c]{r|l|l|l|l|l}
12&\Q(-10)&&\Q(-11)&&\Q(-12)\\\hline
11&&&&&\\\hline
10&&&&&\\\hline
9 &\Q(-8)^2&&\Q(-9)^2&&\Q(-10)^2\\\hline
8 &&&&&\\\hline
7 &&&&&\\\hline
6 &\Q(-4)+\Q(-6)&&\Q(-5)+\Q(-7)&&\Q(-6)+\Q(-8)\\\hline
5 &&&&&\\\hline
4 &&&&&\\\hline
3 &\Q(-2)^2&&\Q(-3)^2&&\Q(-4)^2\\\hline
2 &&&&&\\\hline
1 &&&&&\\\hline
0 &\Q&&\Q(-1)&&\Q(-2)\\\hline
  &0&1&2&3&4
\end{array}$$
\end{table}

By Theorem~\ref{inco}, we know that the rational cohomology of $\mathcal I_1$ is the tensor product of the cohomology of $\GL(3)$ and that of $\mathcal Q_{1}$, hence the spectral sequence given in Table~\ref{tabM31} converges to the tensor product of the cohomology of $\PGL(3)$ and $\mathcal Q_1$. This implies that all non-trivial differentials between column 0 and column 4 in the spectral sequence have rank 1. This yields the following.

\begin{proposition}
The Poincar\'e-Serre polynomial of $\mathcal I_1$ is $(1+\ut2{})(1+\ut43)(1+\ut65)(1+\ut22)(1+\ut{12}6)$.
\end{proposition}

\begin{corollary}\label{M31}
The Poincar\'e-Serre polynomials of $\mathcal Q_{1}$ and of $\Mm31$ are, respectively,
$(1+\ut22)(1+\ut{12}6)$ and $(1+\ut22)(1+\ut22+\ut{12}6)$.
\end{corollary}

\proof
This follows from the above Proposition, together with Theorem~\ref{inco}, and the fact that $\coh\pu{\Hh{3,1}}\cong\coh\pu{\Pp1}$.
\qed

Note that Corollary~\ref{M31} refines the result on the Poincar\'e-Serre polynomial of $\Mm31$ given in \cite{GL}, which corrects an error in \cite[4.10]{Looij}. 

The situation with curves of genus 4 is similar as that in the case of genus 3. The canonical model of a non-hyperelliptic curve of genus 4 is the complete intersection of a quadric and a cubic surface in $\Pp3$ (\cite[IV.5.5.2]{Hartshorne}). As a consequence, we can stratify $\M4$ by taking the space $\M4^0$ of curves which are the complete intersection of a cubic surface and a quadric of maximal rank, the space $\M4^1$ of curves which are the complete intersection of a cubic surface and a quadric cone, and the hyperelliptic locus $\Hh4$. Then, by Remark~\ref{casiLH}, case \ref{4co} and \ref{4ge}, we can determine the cohomology of both $\M4^1$ and $\M4^0$ by computing that of the complements of the discriminants involved in each case.

This construction allows us to compute the rational cohomology of $\M4$. Let us consider $\M4^0$ first. Every element of $\M4^0$ has a canonical embedding to a curve of type $(3,3)$ on $\Pp1\times\Pp1$. The cohomology of $V\setminus\Sigma$ can be obtained by applying Vassiliev-Gorinov's method (as described in Section~\ref{VG}), starting from a classification of all possible singular loci of a $(3,3)$-curve. This yields the following result:

\begin{proposition}[{\cite[1.2]{OT}}]
The cohomology of the moduli space of non-singular $(3,3)$-curves on $\Pp1\times\Pp1$ has Poincar\'e-Serre polynomial equal to $1+\ut65$.
\end{proposition}

Analogously, it is possible to apply Vassiliev-Gorinov's method to the space of non-singular curves of degree 6 on a quadric cone in $\Pp3$ and prove 
that its cohomology is isomorphic to that of the automorphism group of the cone. As a consequence, we have the following.

\begin{theorem}[{\cite[1.4]{OT}}]\label{cohM4}
The Poincar\'e-Serre polynomial of the cohomology of $\M4$ is $1+\ut22+\ut44+\ut65$.
\end{theorem}

\section{Vassiliev-Gorinov's method}\label{VG}

The examples in Section \ref{M3} show that we can determine the rational cohomology of several moduli spaces of smooth curves, once we compute the cohomology of the complement of certain discriminants. There is a specific method, due to Vassiliev, which applies to the latter sort of computation (see e.g. \cite{Vart}). In Section \ref{VGmethod} we present a modification of this method, by Gorinov (see \cite{Gorinov}; for further reference on the construction cited in this section, see the exposition of the method in \cite{OT}). This version of Vassiliev's method applies well to the cases we are interested in. In Section~\ref{qcp} we show how the method works in the special case of the space of quartic plane curves through a fixed point.

In this paper we make an extensive use of Borel-Moore homology, i.e., homology with locally finite support. A reference for its definition and properties is \cite[Chapter 19]{Fulton}.

\subsection{The method}\label{VGmethod}
Here we explain Vassiliev-Gorinov's method for computing the cohomology of complements of discriminants. 

Recall that we are working with a linear space $V$ of sections of a line bundle $L$ on a projective variety $Z$. Our aim is to  compute the rational cohomology of the complement of the discriminant, $X=V\setminus \Sigma$. This is equivalent to determining the Borel-Moore homology of the discriminant, because there is an isomorphism between the reduced cohomology of $X$ and Borel-Moore homology of $\Sigma$. If we denote by $M$ the dimension of $V$, this isomorphism can be formulated as
$$\tilde H^{\pu}(X;\Q)\cong \bar H_{2M-\pu-1}(\Sigma;\Q)(-M).$$

\begin{definition}
A subset $S\subset Z$ is called a \emph{configuration} in $Z$ if it is compact and non-empty. The space of all configurations in $Z$ is denoted by $\Conf(Z)$.
\end{definition}

\begin{proposition}[\cite{Gorinov}]
The Fubini-Study metric on $Z$ induces in a natural way on $\Conf(Z)$ the structure of a compact complete metric space.
\end{proposition}

To every element in $v\in V$, we can associate its singular locus $K_v\in\Conf(Z)\cup\{\emptyset\}$. We have that $K_0=Z$, and that $L(K):=\{v\in V: K\subset K_v\}$ is a linear space for all $K\in \Conf(Z)$. 

Vassiliev-Gorinov's method is based on the choice of a collection of families 
of configurations $X_1,\ldots,X_N\subset\Conf(Z)$, satisfying some axioms (\cite[3.2]{Gorinov}, \cite[List 2.1]{OT}). Intuitively, we have to start by classifying all possible singular loci of elements of $V$. Note that singular loci of the same type will have a space $L(K)$ of the same dimension. We can put all singular configurations of the same type in a family. Then we order all families we get according to the inclusion of configurations. In this way we get a collection of families of configurations which may already satisfy Gorinov's axioms. If this is not the case, the problem can be solved by adding new families to the collection. Typically, the elements of these new families will be degenerations of configurations already considered. 

For instance, configurations with three points on the same projective line and a point outside it can degenerate into configurations with four points on the same line, even if there is no $v\in V$ which is only singular at four collinear points.  

Once the existence of a collection $X_1,\dots,X_N$ satisfying Gorinov's axioms is established, we can proceed with the construction of a resolution for the discriminant. We will formulate this by using the language of \emph{cubical spaces}.

For every $j\in\Z_{\geq1}$, we will denote by $\underline j$ the set $\{1,2,\dots,j\}$. 

\begin{definition}
A \emph{cubical space} over the index set ${\underline{N}}$ (briefly, an ${\underline{N}}$-cubical space) is a collection of topological spaces $\{Y(I)\}_{I\subset{\underline{N}}}$ such that for each inclusion $I\subset J$ we have a natural continuous map $f_{IJ}: Y(J)\rightarrow Y(I)$ such that $f_{IK}=f_{IJ}\circ f_{JK}$ whenever $I\subset J\subset K$. 
\end{definition}

Let us define the cubical spaces we work with.
$$\Lambda(I):=\{A\in\prod_{i\in I} X_i:i<j\Rightarrow A_i\subset A_j\}.$$
$$\X(I):=\{(F,A)\in\Sigma\times\Lambda(I):K_F\supset A_{\max(I)}\}\text{ if }I\neq\emptyset,$$
$$\X(\emptyset):=\Sigma.$$

The subsets $I\subset{\underline{N}}$  are in one to one correspondence with the faces $\Delta_I$ of the simplex
$$\Delta_{\underline{N}}:=\left\{(f:{\underline{N}}\rightarrow[0,1]):\sum_{a\in{\underline{N}}}f(a)=1\right\}.$$
Its faces are defined as $\Delta_I=\{f\in\Delta_{\underline{N}}:f|_{{\underline{N}}-I}=0\}$.  
Note that $\Delta_{\emptyset}=\emptyset$. 
Whenever there is an inclusion $I\subset J$, we can associate to it the map $e_{IJ}:\Delta_I\rightarrow\Delta_J$ given by the inclusion of $\Delta_I$ in $\Delta_J$

\begin{definition}
Let $Y(\pu)$ be a cubical space over the index set ${\underline{N}}$.
Note that $Y(\pu)$ has a natural augmentation towards $Y(\emptyset)$. 
Then the \emph{geometric realization} of $Y(\pu)$ is defined as the map
$$|\epsilon|: \ba Y(\pu)\ba \longrightarrow Y(\emptyset)$$
induced from the natural augmentation on the space
$$\ba Y(\pu)\ba =\bigsqcup_{I\subset{\underline{N}}}(\Delta_I\times Y(I))/R,$$
where $R$ is the equivalence relation generated by 
$$(f,y)\ R\ (f',y')\Leftrightarrow f'=e_{IJ}(f),y=f_{IJ}(y').$$
\end{definition}

Since $\Delta_\emptyset=\emptyset$, the space $Y(\emptyset)$ does not appear  in the construction of the domain of the geometric realization.

For a very natural choice of the topology on $\ba\X(\pu)\ba$ (see \cite{OT}), we have the following:

\begin{proposition}[\cite{Gorinov}]
The geometric realization of $\X(\pu)$,
$$|\epsilon|: \ba \X(\pu)\ba \longrightarrow \X(\emptyset)=\Sigma,$$
is a homotopy equivalence and induces an isomorphism on Borel-Moore homology groups.
\end{proposition}

For every $h$, $1\leq h\leq N$, we have $\underline{h}\hookrightarrow {\underline{N}}$, hence we can restrict $\Lambda(\pu)$ and $\X(\pu)$ to the index set $\underline{h}$, getting the two $\underline{h}$-cubical spaces $\Lambda|_{\underline{h}}(\pu)$ and $\X|_{\underline{h}}(\pu)$. Then for every $1\leq h_1\leq h_2\leq N$, there are natural embeddings $\ba \Lambda|_{\underline{h_1}}(\pu)\ba \hookrightarrow\ba \Lambda|_{\underline{h_2}}(\pu)\ba $ and $\ba \X|_{\underline{h_1}}(\pu)\ba \hookrightarrow\ba \X|_{\underline{h_2}}(\pu)\ba $. 

In this way we can define an increasing filtration on $\ba\Lambda(\pu)\ba $ by posing
$$\Fil_j\ba \Lambda(\pu)\ba :=\im\left(\ba \Lambda|_{\underline{j}}(\pu)\ba \hookrightarrow\ba \Lambda(\pu)\ba \right)$$
for $j=1,\dots,N$. We define analogously the filtration $\Fil_j\ba \X(\pu)\ba $ on $\ba \X(\pu)\ba $. We use the notation $F_j:=\Fil_j\ba \X(\pu)\ba\setminus\Fil_{j-1}\ba \X(\pu)\ba $ and
$\Phi_j:=\Fil_j\ba \Lambda(\pu)\ba\setminus\Fil_{j-1}\ba \Lambda(\pu)\ba $ for these strata.

The filtration $\Fil_j\ba \X(\pu)\ba $ defines a spectral sequence that converges to the Borel-Moore homology of $\Sigma$. Its term $E^1_{p,q}$ is isomorphic to $\bar H_{p+q}(F_p;\Q)$.

\begin{proposition}[\cite{Gorinov}]\label{ucci}
\renewcommand{\labelenumi}{\arabic{enumi}.}
\begin{enumerate}
\item For every $j=1,\dots,N$, the stratum $F_j$ is a complex vector bundle  over $\Phi_j$. The space $\Phi_j$ is in turn a fiber bundle over the configuration space $X_j$. 
\item If $X_j$ consists of configurations of $m$ points, the fiber of $\Phi_j$ over any $x\in X_j$ is an $(m-1)$-dimensional open simplex, which changes its orientation under the homotopy class of a loop in $X_j$ interchanging a pair of points in $x_j$.
\item\label{opeco}  If $X_N=\{Z\}$, $F_N$ is the open cone with vertex a point (corresponding to the configuration $Z$),  over $\Fil_{N-1}\ba \Lambda(\pu)\ba $. 
\end{enumerate}\end{proposition}

We recall here the topological definition of an open cone.
\begin{definition}
Let $B$ be a topological space. Then a space is said to be an \emph{open cone} over $B$ with vertex a point if it is homeomorphic to the space $B\times[0,1)/R$, where the equivalence relation is $R=(B\times\{0\})^2$.
\end{definition}

The fiber bundle $\Phi_j\rightarrow X_j$ of Proposition \ref{ucci} is in general non-orientable. As a consequence, we have to consider the homology of $X_j$ with coefficients not in $\Q$, but in some local system of rank one. Therefore we recall some constructions concerning Borel-Moore homology of configuration spaces with twisted coefficients.

\begin{definition}
Let $Z$ be a topological space. Then for every $k\geq 1$ we have the space of ordered configurations of $k$ points in $Z$,
$$\F kZ=Z^k\setminus\bigcup_{1\leq i<j\leq k}\{(z_1,\dots,z_k)\in Z^k: z_i=z_j\}.$$
There is a natural action of the symmetric group $\s_k$ on $F(k,Z)$. The quotient is called the space of unordered configurations of $k$ points in $Z$,
$$\B kZ=\F kZ/\s_k.$$
\end{definition}

The \emph{sign representation} $\pi_1(\B kZ)\rightarrow \Aut(\Z)$ maps the paths in $\B kZ$ defining odd (respectively, even) permutations of $k$ points to multiplication by $-1$ (respectively, 1). The local system $\pm\Q$ over $\B kZ$ is the one locally isomorphic to $\Q$, but with monodromy representation equal to the sign representation of $\pi_1(\B kZ)$. We will often call $\bar H_\pu (\B kZ,\pm\Q)$ the \emph{Borel-Moore homology of $\B kZ$ with twisted coefficients}, or, simply, the \emph{twisted Borel-Moore homology of $\B kZ$}.

\subsection{Quartic curves through a point}\label{qcp}
In this section, we will give an example of how to apply Vassiliev-Gorinov's method to a concrete case, namely, to the case of quartic curves through a fixed point. This is an adaptation of the computation of Vassiliev in \cite{Vart} for the space of quartic curves.

Let us fix a point $p\in\Pp2$, and choose coordinates $x_0,x_1,x_2$ on $\Pp2$. Denote by $V$ the vector subspace of $\C[x_0,x_1,x_2]_4$ of homogeneous polynomials of degree four vanishing at $p$. The singular polynomials inside $V$ form the discriminant locus $\Sigma\subset V$.  We want to apply Vassiliev-Gorinov's method in order to determine the Borel-Moore homology of $\Sigma$ and the cohomology of its complement in $V$.

We start by listing in Table~\ref{l1} all possible singular loci of elements of $\Sigma$. Note that the list of configurations in Table~\ref{l1} defines a collection of families of configurations that already satisfies Gorinov's axioms. Therefore, we can use this collection to construct the cubical spaces $\X(\pu)$ and $\Lambda(\pu)$ as in Section~\ref{VGmethod}.

\begin{table}[t]\caption{Singular sets of quartic polynomials vanishing at $p$ \label{l1}.}
The number in the last column is the dimension of the linear subspace of $V$ of polynomials singular in a configuration of the given type. The position of $n$ points is called \emph{general} when no three of them lie on the same line. 

\begin{center}
\begin{tabular}[c]{|p{0.9cm}|p{8.8cm}|l|} 
\hline
Type&Description&Dim.\\\hline
{1a}& The point $p$. & 12\\\hline
{1b}& Any other point. & 11\\\hline
{2a}& The point $p$ and another point. & 9\\\hline
{2b}& Any other pair of points. & 8\\\hline
{3a}& Three points lying on a line through $p$. &7\\\hline
{3b}& Three points on a line not passing through $p$. &6\\\hline
{4a}& The point $p$ and two other points in general position with respect to $p$.& 6 \\\hline
{4b}& Three points in general position, different from $p$. & 5\\\hline
{5a}& A line through $p$.& 6\\\hline
{5b}& A line not passing through $p$.& 5\\\hline
{6a,b,c}& Three points on the same line and one point outside it. There are 3 subcases, according to the position of $p$.& 4,4,3\\\hline
{7a,b}& Four points in general position. There are 2 subcases, depending on the position of $p$.& 3,2\\\hline
{8a,b,c}& A line and a point outside it. There are 3 subcases, depending on the position of $p$.& 3,3,2
\\\hline
{9a,b,c}& Two pairs of points, lying on different lines, and the intersection point of the two lines they span. There are 3 subcases: 
$p$ may be the central point, or lie on one of the lines, or the conic is uniquely determined as the conic passing through $p$ and the other 4 points.& 2,2,1 \\\hline
{10}& Six points, pairwise intersection of a general configuration of four lines in $\Pp2$. At least one line passes through $p$.& 1\\\hline
{11}& A non-singular conic through $p$.& 1\\\hline
{12}& A conic of rank 2, passing through $p$.&1 \\\hline
{13}&The whole $\Pp2$.& 0\\\hline
\end{tabular}
\end{center}
\end{table}

The strata defined by configurations of types 3a, 3b, 4b, 5, 6, 7, 8, 9, 11 and 12 do not contribute to the Borel-Moore homology of the domain of the geometric resolution of $\X(\pu)$. This follows either from the fact that the twisted Borel-Moore homology of the corresponding configuration spaces is trivial, or from the fact that the configurations are curves. It remains to calculate the Borel-Moore homology of the other strata.

\begin{itemize}
\item[1a] The space $F_{1a}$ is isomorphic to $\C^{12}$. We have $\Phi_{1a}\cong \{p\}$.
\item[1b] The space $F_{1b}$ is a complex vector bundle of rank 11 over $\Phi_{1b}\cong(\Pp2\setminus\{p\})$.
\item[2a] The space $F_{2a}$ is a complex vector bundle of rank 9 over $\Phi_{2a}$, which in turn is a non-orientable $\Delta_1$-bundle over $\Pp2\setminus\{p\}$.
\item[2b] The space $F_{2b}$ is a complex vector bundle of rank 8 over $\Phi_{2b}$, which in turn is a non-orientable $\op\Delta_1$-bundle over $\B2{\Pp2\setminus\{p\}}$. After the choice of a line through $p$, the complement of $p$ in $\Pp2$ can be seen as the disjoint union of an affine line and an affine two-space, hence we have $\B2{\Pp2\setminus\{p\}}=\B2{\C}\sqcup \C\times\C^2\sqcup \B2{\C^2}$. The induced local system is $\pm\Q$ on the first and the third component. It is a standard fact in the theory of complements of unions of hyperplanes that $\bar H_\pu(\B j{\C^N};\Q)=0$ when $j\geq2$. Hence,
$\bar H_\pu(\B2{\Pp2\setminus\{p\}};\pm\Q)$ $\cong\bar H_\pu(\C^3;\Q)$. 
\item[4a] The space $F_{4a}$ is a complex vector bundle of rank 6 over $\Phi_{4a}$, which in turn is a non-orientable $\op\Delta_2$-bundle over $\B2{\Pp2\setminus\{p\}}$.
\item[10] The space $F_{10}$ is a complex vector bundle of rank 1 over $\Phi_{10}$, which in turn is a non-orientable $\op\Delta_5$-bundle over the configuration space $X_{10}$. 
 
It is easier to consider $X_{10}$ as a subspace of $\B4{\Pp2\duale}$. Then $X_{10}$ is the space of all configurations of four points in general position such that at least one point lies on the line $p\duale$. The local system induced by the orientation of the simplicial bundle is the constant one, because interchanging two lines interchanges two pairs of intersection points. 

The Borel-Moore homology with constant coefficients of $\tB4{\Pp2\duale}$ is easy to compute, because it coincides with that of $\PGL(3)$. Its Poincar\'e-Serre polynomial is $\ut{-16}{16}+\ut{-12}{13}+\ut{-10}{11}+\ut{-6}{8}$.

Consider the following variety:
$$\mathcal S:=\{(l_1,l_2,l_3,l_4,p)\in\tF4{\Pp2}\times\Pp2: p\in\bigcup_il_i\},$$
where $\tF4{\Pp2\duale}$ denotes the space of four lines in $\Pp2$  in general position.
If we consider the projection of $\mathcal S$ onto $\tF4{\Pp2\duale}$, it is easy to see that the $\s_4$-invariant part of the Borel-Moore homology of $\mathcal S$ coincides with the tensor product of the Borel-Moore homology of $\tF4{\Pp2\duale}\cong\PGL(3)$ and that of $\Pp1$. If we consider the projection $p_2:\mathcal S\longrightarrow\Pp2$, we obtain that the fibre of $p_2$ is a $\s_4$-covering of $X_{10}$. This implies that the cohomology of $X_{10}$ coincides with the $\s_4$-invariant part of the cohomology of this fibre. A direct computation based on the Leray spectral sequence associated to $p_2$ yields that the Poincar\'e-Serre polynomial of the Borel-Moore homology of $X_{10}$ must be $\ut{-14}{14}+2\ut{-10}{11}+\ut{-6}8.$

\item[13] The space $F_{13}$ is an open cone with base $\Fil_{13}\ba\Lambda(\pu)\ba=\bigcup_{1\leq i\leq13}\Phi_i$. The spectral sequence associated to the filtration $\Fil_i\left(\ba\Lambda(\pu)\ba\right)$ is such that the differential $d_1$ is an isomorphism between the columns corresponding to the contributions of $\Phi_{2a}$ and $\Phi_{1b}$. The same holds for $\Phi_{4a}$ and $\Phi_{2b}$.
 From this it follows that the only contribution comes from configurations of type 10, and that $\bar H_\pu(F_{13};\Q)\cong \bar H_{\pu-1}(\Phi_{10};\Q)$.
\end{itemize}

The spectral sequence converging to the Borel-Moore homology of $\Sigma$ has $E^1$ term as in Table~\ref{tabsigma}.

\begin{table}[t]\caption{\label{tabsigma}Spectral sequence converging to $\bar H_\pu(\Sigma;\Q)$.}
$$\begin{array}[c]{r|l|l|l|l|l|l|l}
24&      &\Q(13)&      &      &     &&\\\hline
23&\Q(12)&      &      &      &     &&\\\hline
22&      &\Q(12)&      &      &     &&\\\hline
21&      &      &      &      &     &&\\\hline
20&      &      &\Q(11)&      &     &&\\\hline
19&      &      &      &\Q(11)&     &&\\\hline
18&      &      &\Q(10)&      &     &&\\\hline
17&      &      &      &      &     &&\\\hline
16&      &      &      &      &     &&\\\hline
15&      &      &      &      &\Q(9)&\Q(8)&\\\hline
14&      &      &      &      &     &&\\\hline
13&      &      &      &      &     &&\Q(7)\\\hline
12&      &      &      &      &     &\Q(6)^2&\\\hline
11&      &      &      &      &     &&\\\hline
10&      &      &      &      &     &&\Q(5)^2\\\hline
 9&      &      &      &      &     &\Q(4)&\\\hline
 8&      &      &      &      &     &&\\\hline
 7&      &      &      &      &     &&\Q(3)\\\hline
 6&      &      &      &      &     &&\\\hline
  &1a    &1b    &2a    &2b    &4a   &10&13
\end{array}$$
\end{table}

From the above, and considering the fact that the relationship between Borel-Moore homology of $\Sigma$ and reduced cohomology of $X=V\setminus\Sigma$ is given by 
$$\tilde H^{\pu}(X;\Q)\cong \bar H_{27-\pu}(\Sigma;\Q)(-14),$$
we get the Poincar\'e-Serre polynomial of the
rational cohomology of $V\setminus\Sigma$
\begin{equation}\label{I1}
(1+\ut2{})(1+2\ut43+\ut86+\ut{12}6+2\ut{16}9+\ut{20}{12}).
\end{equation}

This result is used in Section~\ref{M3} to compute the cohomology of $\Mm 31$.

\end{document}